%% file: LD-GTransNet_rev.tex
\documentclass[preprint,11pt]{elsarticle}
\usepackage{amssymb,amsfonts,amsthm,bm}
\usepackage{amsmath}
\usepackage{mathtools}
\usepackage{algorithm}
\usepackage{algpseudocode}
\usepackage{booktabs}
\usepackage{multirow}
\usepackage{colortbl}
\usepackage{diagbox}
\usepackage{graphicx}
\usepackage[english]{babel}
\usepackage{listings}
\usepackage[colorlinks=true, allcolors=blue]{hyperref}
\usepackage{enumerate}
\usepackage{makecell}

\usepackage{subcaption}
\usepackage{arydshln}
\usepackage{xargs}
\usepackage{float}
\usepackage{cases}
\usepackage{float}
\usepackage{fullpage}

\theoremstyle{plain}
 
\newtheorem{theorem}{Theorem}[section] 
 
\newtheorem{remark}{Remark}[section] 

\input{preamble}

\graphicspath{{./figure/}}

\begin{document}

\begin{frontmatter}
    \title{Efficient and Accurate Lift-and-Decoupling Generalized Transferable Neural Networks  for Elliptic and Moving Interface Problems}

    \author[whu]{Tao Cheng} 
    \ead{taocheng77@whu.edu.cn}
    
    \author[whu]{Junqing Jia} 
    \ead{whujjq@whu.edu.cn}
    
    \author[usc]{Lili Ju\corref{cor1}} 
    \ead{ju@math.sc.edu}
    \cortext[cor1]{Corresponding author}
    
    \author[whu]{Xiaoping Zhang} 
    \ead{xpzhang.math@whu.edu.cn}

    \affiliation[whu]{organization={School of Mathematics and Statistics, Wuhan University},
        city={Wuhan},
        postcode={430072},
        state={Hubei},
        country={China}}

    \affiliation[usc]{
      organization={Department of Mathematics, University of South Carolina}, 
      city={Columbia},
      postcode={29208}, 
      state={SC},
      country={USA}}

    \begin{abstract}
        In this paper, we propose a highly efficient and accurate method, called ``lift-and-decoupling generalized transferable neural network" (LD-GTransNet), for numerical solution of elliptic and moving interface problems. Our method makes use of the multi-layer GTransNet approach and incorporates a novel {lift-and-decoupling} strategy, in which the hidden-layer neuron parameters are preset and the physical coordinates and the added auxiliary variable are processed by separate subnetworks. Such design enables a unified global approximation of piecewise smooth solutions without explicit domain decomposition, while preserving network expressivity and avoiding the redundancy inherent in conventional {lift} approaches. The resulting formulation reduces to a linear least-squares problem, thereby eliminating the need of nonlinear training and ensuring the high computational efficiency and accuracy.  Extensive numerical experiments on benchmark interface problems in two and three dimensions demonstrate that, compared with existing neural network-based and traditional numerical methods, our proposed LD-GTransNet  consistently delivers superior performance, particularly for problems with high-contrast coefficients and complex interface geometries.
    \end{abstract}




\end{frontmatter}

\section{Introduction}\label{sec:1}
Interface problems arise in a wide range of scientific and engineering applications involving heterogeneous media and phase transitions. Typical examples include heat conduction in composite materials, multiphase flows with surface tension effects, electrodiffusion across biological membranes, and electromagnetic wave propagation in discontinuous media \citep{Choi2021,Gibou2019,Kharangate2017}. These problems are characterized by domains partitioned into subregions by lower-dimensional interfaces, across which the solution and/or its flux exhibit discontinuities governed by prescribed jump conditions. Such discontinuities often induce low regularity or even singular behavior in the solution, posing significant challenges for numerical approximation.


Classical numerical methods for interface problems can be broadly categorized into body-fitted and unfitted approaches. Body-fitted methods align the computational mesh with the interfaces, allowing discontinuities to be explicitly resolved within the discretization framework. Representative examples include finite element \citep{Babu1970, Bramble1996}, finite volume \citep{Guittet2015}, virtual element \citep{Chen2017}, and discontinuous Galerkin methods \citep{Cai2011}. While these approaches can achieve high-order accuracy when the mesh conforms well to the interface geometry, generating high-quality fitted meshes remains a major bottleneck, especially for problems with complex geometries, moving interfaces, or multiscale features.
To alleviate the dependence on mesh generation, unfitted methods have been extensively developed. These methods employ interface-independent meshes and incorporate jump conditions through modified discretizations, enriched approximation spaces, or penalty formulations. Notable examples include immersed interface methods \citep{LeVeque1994}, immersed finite element methods \citep{Lin2015, Guo2019, Guo2020}, and unfitted discontinuous Galerkin and finite volume methods \citep{Chen2008, Coco2018, Massjung2012, Thacher2023,Gibou2020}. Despite their flexibility, unfitted methods often require sophisticated treatments near cut cells and accurate geometric reconstruction, which can become highly challenging for interfaces with large curvature or intricate topology, potentially leading to reduced robustness and accuracy.


The fundamental limitation of both fitted and unfitted approaches lies in their reliance on mesh-based representations of geometry. This limitation has motivated the development of meshfree methods, which approximate solutions directly on scattered nodes without explicit mesh construction. Among them, kernel-based methods \citep{Kraus2023} and generalized finite difference schemes \citep{Davydov2021, Martin2017_2, Martin2017_1} have shown promising performance. More recently, neural network-based approaches \citep{Wang2020} have emerged as a powerful alternative due to their strong approximation capability and inherent flexibility in handling high-dimensional and irregular problems.


Neural network methods for PDEs can be broadly divided into two categories according to their training paradigms. The first category consists of fully parameterized networks, such as the Deep Galerkin Method \citep{Sirignano2018}, Deep Ritz Method \citep{Yu2018}, and Physics-Informed Neural Networks (PINNs) \citep{Raissi2019}, where all network parameters are optimized via gradient-based algorithms. Although these methods are highly expressive, they often suffer from expensive training, sensitivity to hyperparameters, and difficulties in optimization due to non-convexity. Various strategies, including domain decomposition \citep{Wu2022, Roy2025, Zhu2023}, weak formulations \citep{Wang2020}, and input augmentation and embedding \citep{Kumar2025, Hu2022, Tseng2023}, have been proposed to address interface problems within this framework. However, these approaches typically introduce additional complexity in network design or training. The second category, referred to as fixed-neuron methods, offers a computationally efficient alternative. In these approaches, the parameters of hidden layers are  preset and fixed, while only the output weights are determined through the least-squares optimization. Representative methods include the Extreme Learning Machine (ELM) \citep{Calabr2021, Dong2021, Quan2023, Zeng2025}, random feature methods (RFM) \citep{Chen2023, Chi2024, Sun2025}, and transferable neural networks (TransNet) \citep{Lu2025, Shen2026, Zhang2024,Sun2026,cheng2026BDFTC}. By avoiding or reducing nonlinear training, these methods significantly improve efficiency and robustness. Nevertheless, their direct application to interface problems still remains challenging due to two key issues: (i) the lack of sufficiently expressive basis functions for capturing sharp discontinuities, and (ii) the difficulty of incorporating interface conditions without reintroducing geometric complexity.

Existing neural network approaches for interface problems can be broadly categorized into two types. 
{\em Domain decomposition methods} partition the computational domain along physical interfaces and assign separate networks to each subdomain, enforcing jump conditions via additional loss terms \citep{Chi2024,Lu2025,Wu2022}. While conceptually straightforward, these methods suffer from increasing architectural complexity and computational cost as the number of subdomains grows, and they require careful balancing of multiple loss components.
{{\em Dimension lift methods} offer an alternative route by introducing an auxiliary variable to lift a low-dimensional piecewise smooth solution into a higher-dimensional function space, thereby transforming the problem into one with global smoothness.}  Specifically, an auxiliary variable \(\xi = \xi(\bx)\) acts as a regional indicator, assigning distinct constant values to each subdomain. A globally smooth function \(F(\bx, \xi)\) is then constructed such that, by fixing \(\xi\) to a subdomain label, the physical solution is recovered \citep{Hu2022,Lai2026,SONG2026,Lai2026}. This approach avoids explicit domain decomposition and enables a single global network to approximate the solution. {However, directly concatenating the physical coordinates and the auxiliary variable as inputs to a single network leads to an unnecessary increase in input dimensionality. Moreover, different \((\bx, \xi)\) combinations may correspond to the same physical state, giving rise to representation non-uniqueness and optimization difficulties.}


In this paper, we propose a lift-and-decoupling generalized transferable neural network (LD-GTransNet) to overcome the above limitations. Our approach builds upon the GTransNet \cite{cheng2026GTransNet}, a multi-layer extension of TransNet that generates structured and expressive neural basis functions through a theoretically guided. To effectively handle discontinuities, we introduce a lifting framework with a key modification: {unlike existing approaches that directly lift all variables into a unified high-dimensional space,}  we decouple the physical coordinates and the auxiliary variable into two parallel subnetworks. Such design preserves the core advantages of lifting while removing unnecessary increases in model complexity and mitigating the issues caused by representation non-uniqueness. We particularly emphasize that, compared with existing dimension lift methods such as DCSNN~\cite{Hu2022} and DC-RFM~\cite{SONG2026}, the proposed LD-GTransNet exhibits an essential distinction: existing methods typically concatenate auxiliary variables with spatial inputs and embed them into a unified high-dimensional space, thereby producing a large number of redundant basis functions that are never activated simultaneously; by contrast, our decoupled architecture processes spatial coordinates and auxiliary labels through parallel subnetworks, and the Hadamard product of their outputs selectively activates spatial basis functions according to the subdomain label. This mechanism effectively eliminates representational redundancy while maintaining full approximation capability across all subdomains. {Furthermore, when the original problem contains no interface (i.e., the solution is globally smooth), the auxiliary variable degenerates to a constant, and our method naturally reduces to the standard GTransNet, exhibiting favorable backward compatibility.}


The remainder of this paper is organized as follows. Section \ref{sec:2} reviews the GTransNet architecture and its use  for solving PDEs. Section \ref{sec:3} presents in detail the proposed LD-GTransNet method for elliptic and moving interface problems. Section \ref{sec:4} provides extensive numerical experiments and comparisons to demonstrate superior accuracy and efficiency of our LD-GTransNet. Finally, some conclusions and future work are discussed in Section \ref{sec:5}.


\section{Generalized Transferable Neural Networks for PDEs}\label{sec:2}

\subsection{Transferable Neural networks}
The TransNet developed in \cite{Zhang2024} is a two-layer (single hidden layer) neural network method for solving PDEs of the following form
\begin{equation}\label{modelp}
\begin{cases}
\mathcal{L}(u)(\bx) = f(\bx), & \bx \in \Omega, \\
\mathcal{B}(u)(\bx) = g(\bx), & \bx \in \partial\Omega,
\end{cases}
\end{equation}
where $\Omega \subset \mathbb{R}^d$ is an open bounded spatial domain with a Lipschitz boundary \(\partial \Omega\),
$u(\bx):\Omega\rightarrow\mathbb{R}$ denotes the unknown function, \(\mathcal{L}(\cdot)\) defines a differential operator, and \(\mathcal{B}(\cdot)\) specifies 
certain boundary conditions. The key ingredient is the so-called neural feature space, denoted by $\mathcal{P}_{\rm NN}$:
\[
\mathcal{P}_{\rm NN} = {\rm span}\{\psi_1(\bx), \psi_2(\bx), \ldots, \psi_N(\bx)\}
\]
where $\psi_i(\bx) = \sigma(\bw_i^\top \bx + b_i)$ for $i = 1,2,\cdots,N$ are the neural basis functions, where $\sigma(\cdot)$ is a nonlinear activation function, (e.g., $\tanh$) and $\bw_i \in \mathbb{R}^d$, $b_i \in \mathbb{R}$ are the weight and bias of the $i$-th hidden-layer neuron, respectively. The TransNet
solution of \eqref{modelp} $u_{\rm NN} \in \mathcal{P}_{\rm NN}$ is represented by a linear combination of the neural basis functions:
\begin{equation}
u_{\rm NN} = \sum_{i=1}^N \alpha_i \psi_i(\bx).
\end{equation}

The core innovation of TransNet is that its single hidden-layer parameters are pre-determined to achieve a uniform spatial distribution, independent of the specific PDE. To accomplish this, TransNet reparameterizes each neuron as:
\begin{equation}\label{eq:TransNet_Init}
\psi_i(\bx) = \sigma\big(\gamma_i (\mathbf{a}_i^\top \bx + r_i)\big),\quad i=1,\ldots,N,
\end{equation}
Here, the unit vector $\mathbf{a}_i = \mathbf{w}_i/\gamma_i$ (i.e., $\|\mathbf{a}_i\|_2 = 1$, where $\|\cdot\|_2$ denotes the Euclidean norm) and the scalar $r_i = b_i/\gamma_i$ serve as the location parameter defining the direction and the offset of the $i$-th partition hyperplane, respectively. The parameter $\gamma_i > 0$ is a predefined shape factor controlling the steepness of the activation.
Taking the unit ball domain $B_1(\mathbf{0}) \subset \mathbb{R}^d$ as an example, TransNet achieves uniform expressivity by sampling the normal vectors $\mathbf{a}_i$ uniformly from the unit sphere $\mathbb{S}^{d-1}$ (via normalized standard Gaussian vectors):
\begin{equation}\label{eq:sample_a}
\mathbf{a}_i = \frac{\mathbf{z}}{|\mathbf{z}|}, \quad \mathbf{z} \sim \mathcal{N}(\mathbf{0}, {\bf I}_d),
\end{equation}
and the offsets $r_i$ uniformly from $\mathcal{U}[0,1]$. The following theorem formalizes this uniform distribution:

\begin{theorem}[Uniform neuron distribution in the unit ball $B_1(\mathbf{0})$ \cite{Zhang2024}]\label{thm:uniform_neuron}
Given a set of $N$ hyperplanes of $\mathbb{R}^d$ defined by
\[
\mathbf{a}_i^\top \mathbf{x} + r_i = 0, \quad i = 1,2,\dots,N.
\]
If $\{\mathbf{a}_i\}_{i=1}^N$ are i.i.d. and uniformly distributed on the $d$-dimensional unit sphere, and $\{r_i\}_{i=1}^N$ are i.i.d. and uniformly distributed in $[0,1]$, then for a fixed $\tau \in (0,1)$,
\begin{equation}
\mathbb{E}\left[ D_N^\tau(\mathbf{x}) \right] = \tau, \quad \forall \|\mathbf{x}\|_2 \leq 1 - \tau,
\end{equation}
where $D_N^\tau(\mathbf{x})$ is the density function of the neurons defined by
\begin{equation}
D_N^\tau(\mathbf{x}) = \frac{1}{N} \sum_{i=1}^N \chi_{\{d_i(\mathbf{x}) < \tau\}}(\mathbf{x}),
\end{equation}
with $\chi_{\{d_i(\mathbf{x}) < \tau\}}(\mathbf{x})$ denoting the indicator function whether the distance between $\mathbf{x}$ and the $i$-th partition hyperplane is smaller than $\tau$.
\end{theorem}

By Theorem \ref{thm:uniform_neuron}, the hidden-layer neurons are uniformly distributed within $B_1(\mathbf{0})$. For a general problem domain $\Omega$, TransNet utilizes a bounding ball $B_R(\mathbf{x_c})$ slightly larger than $\Omega$, centered at its approximate geometric center $\bx_c$. Through an affine transformation, the corresponding partition hyperplanes are mapped to:
\begin{equation}
\mathbf{a}_i^\top (\bx - \bx_c) + R r_i = 0.
\end{equation}
This ensures uniform expressivity across $\Omega$.

\subsection{Generalized transferable neural networks}

To enhance the applicability and robustness of TransNet in solving PDEs, GTransNet was proposed as a multi-hidden-layer generalization \cite{cheng2026GTransNet}. Unlike the original shallow network, GTransNet employs a deep fully connected architecture with $L \ge 2$ hidden layers. The first hidden layer largely inherits the predetermined uniform sampling scheme of TransNet, with the modification that the neuron biases are drawn from a centrally symmetric distribution. For the subsequent hidden layers, bias terms are entirely omitted, and the weight matrices are predetermined using a variance-controlled sampling strategy. Finally, as in TransNet, the output layer weights are optimized by minimizing a residual-based loss function. 

\begin{figure}[ht!]
    \centering
    \includegraphics[width=0.99\linewidth]{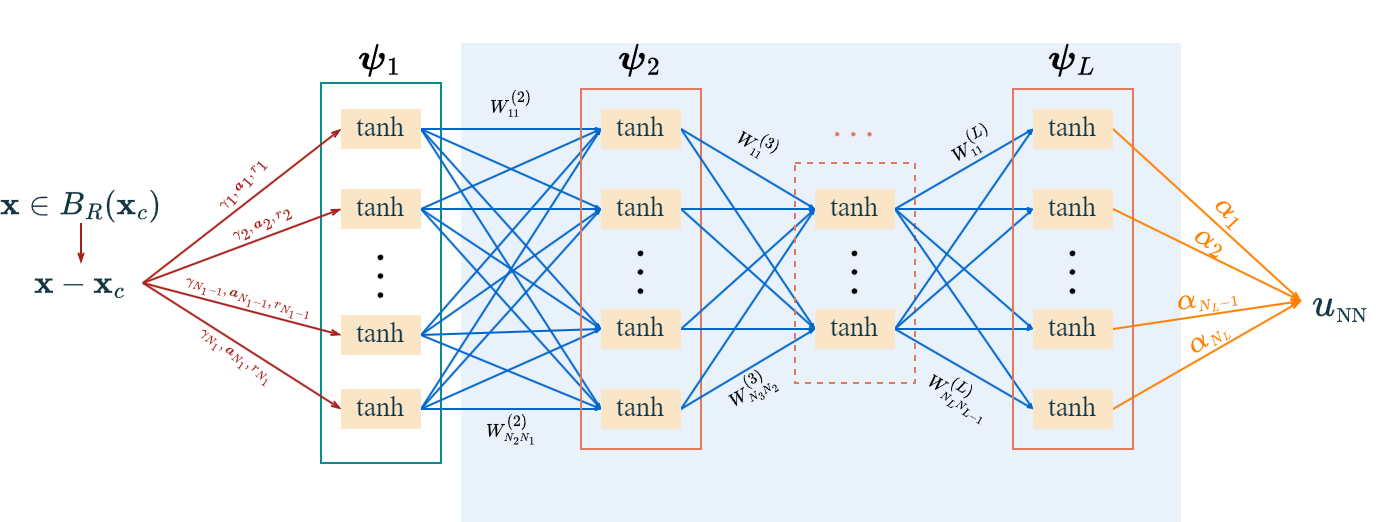}
     \vspace{0.1cm}
    \caption{Network architecture of the proposed GTransNet \eqref{eq:g_transnet} with $L$ hidden layers.}
    \label{fig:GTransNet}
\end{figure}

Assuming the physical domain is bounded by a ball $\Omega \subset B_R(\mathbf{x}_c)$, the forward pass of GTransNet with $L$ hidden layers is formulated as follows:
\begin{equation}\label{eq:g_transnet}
    \begin{cases}
        \vpsi_1 = \tanh(\boldsymbol{\gamma}\odot(\bA(\mathbf{x}-\mathbf{x}_c) + R\br)), \\
        \vpsi_2 = \tanh(\bW_2\vpsi_1), \\
        \vdots\\
        \vpsi_L = \tanh(\bW_L\vpsi_{L-1}),
    \end{cases}
\end{equation}
where 
$\vpsi_1=(\psi^{(1)}_1,\psi^{(1)}_2,\cdots,\psi^{(1)}_{N_1})^\top\in \R^{N_1}$, $\bA = (\ba_1, \ba_2,\ldots, \ba_{N_1})^\top\in\R^{{N_1}\times d}$, $\boldsymbol{\gamma}=(\gamma_1,\gamma_2,\cdots,\gamma_{N_1})^\top\in \R^{N_1}$,  $\br = (r_1, r_2\ldots, r_{N_1})^\top\in\R^{{N_1}}$, $\odot$ denotes the elementwise multiplication, and for $l=2, 3,\ldots, L$, 
$$
\vpsi_l=(\psi^{(l)}_1,\cdots,\psi^{(l)}_{N_l})\in \R^{N_l}, \quad
\bW_l=(W^{(l)}_{ij})\in \R^{N_l\times N_{l-1}}, 
$$ 
where $N_l$ denotes the number of neurons in the $l$-th hidden layer. For the subsequent hidden layers ($l \ge 2$), $\mathbf{W}_l \in \mathbb{R}^{N_l \times N_{l-1}}$ is the predetermined weight matrix without any bias vector. The architecture of the proposed GTransNet is illustrated in Figure \ref{fig:GTransNet}. {For clarity, $\vpsi_l(\bx) \in \R^{N_l}$ denotes the feature vector output by the $l$-th hidden layer, while its individual components $\{\psi_i^{(l)}(\bx)\}_{i=1}^{N_l}$ constitute the set of neural basis functions at layer $l$. In particular, the final-layer features $\{\psi_i^{(L)}(\bx)\}_{i=1}^{N_L}$ form the basis of GTransNet for the solution approximation.}  If $L=1$, then the GTransNet clearly reduces to the original TransNet. 

The  variance-controlled sampling strategy for $\{\bW_l\}_{l=2}^{L}$ is inspired by Xavier initialization \cite{Xavier}, which aims to maintain a consistent signal variance across layers during both forward and backward propagation. Specifically, the entries $W_{ij}^{(l)}$ of $\mathbf{W}_l$ are independently drawn from a zero-mean Gaussian distribution:
\begin{equation}\label{eq:sample_W}
    W_{ij}^{(l)} \sim \mathcal{N} \left(0,\sigma_l^2\right), \quad \text{with} \quad \sigma_l = \sqrt{\frac{\delta}{N_{l-1}}}, \quad l = 2, \dots, L,
\end{equation}
where $N_{l-1}$ is the input dimension of the $l$-th layer, and $0 < \delta \le 1$ serves as a user-defined variance control parameter to adjust the steepness and scaling of the forward signals.
The final network output $u_{\mathrm{NN}}(\mathbf{x})$ is a linear combination of the neural basis functions, i.e. $L$-th layer features,
\begin{equation}
u_{\rm NN} = \sum_{i=1}^{N_L} \alpha_i \psi^{(L)}_i,
\end{equation}
 parameterized by the trainable coefficients $\boldsymbol{\alpha} \in \mathbb{R}^{N_L}$.

\section{Lift-and-Decoupling GTransNets for Elliptic and Moving Interface Problems}\label{sec:3}
The GTransNet framework introduced in Section \ref{sec:2} provides an efficient architecture for solving various partial differential equations with globally smooth solutions, but cannot directly capture solution discontinuities across interfaces. To address this limitation, we propose a highly accurate and efficient  lift-and-decoupling GTransNet, which lifts the piecewise smooth solution into a higher-dimensional space by introducing an auxiliary variable and then adopt a decoupling strategy in the representation of the solution and the construction of the whole network architecture. Section \ref{sec:3.1} presents a detailed description of the 
LD-GTransNet for solving elliptic interface problems, and Section \ref{sec:3.2} extends it to moving interface problems.

\subsection{Elliptic interface problems}\label{sec:3.1}

\begin{figure}[ht!]
    \centering
    \includegraphics[width=0.7\linewidth]{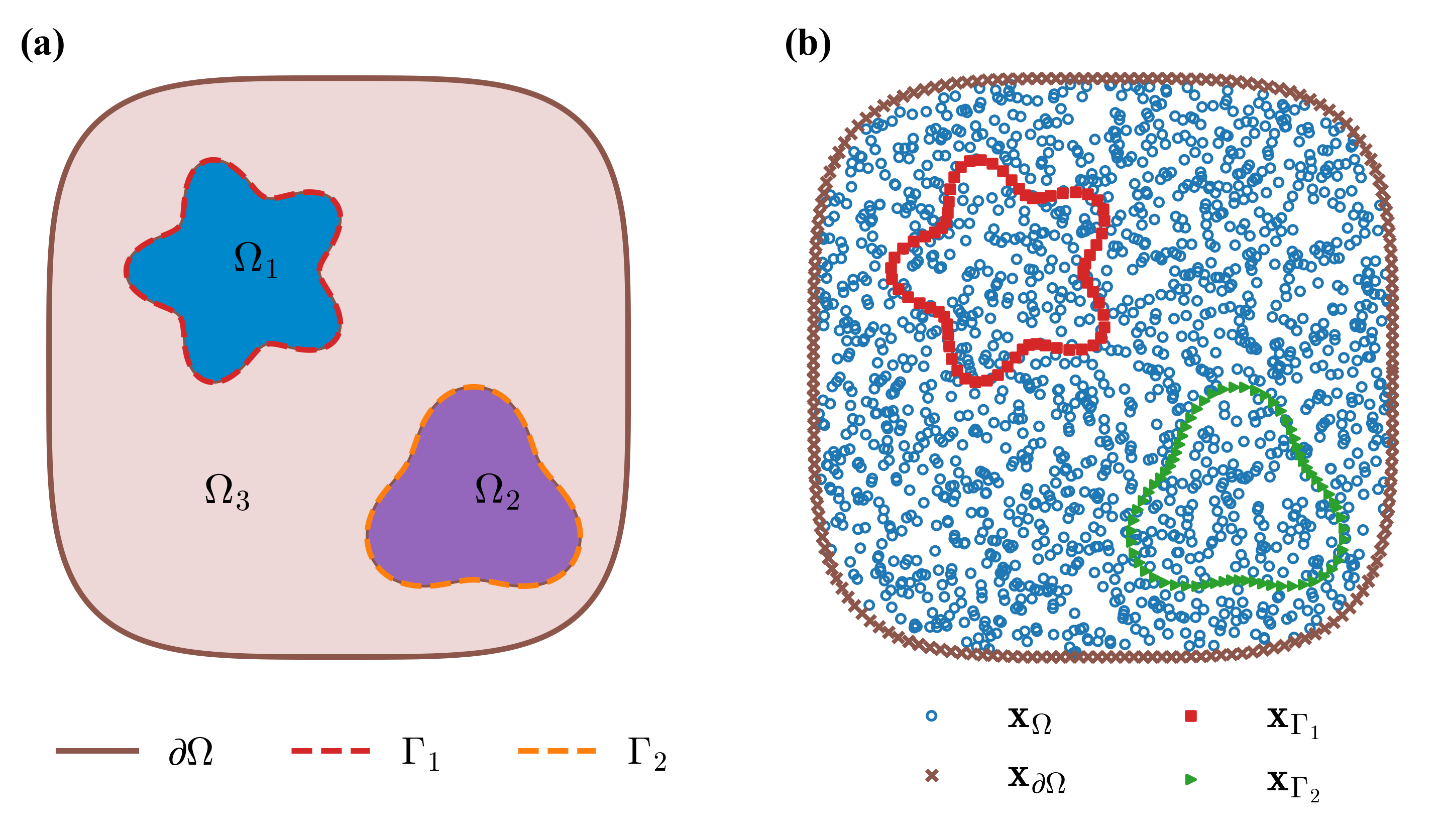}
    \vspace{-0.3cm}
    \caption{(a): domain with interfaces, $\Gamma_1 = \overline{\Omega}_1\cap\overline{\Omega}_3,\Gamma_2 = \overline{\Omega}_2\cap\overline{\Omega}_3$.
    (b): collocation points on the $\Omega,\partial\Omega,\Gamma_1$ and $\Gamma_2$.}
    \label{fig:interface_sketch}
\end{figure}

In this subsection, we consider the following elliptic interface problem defined on an open bounded domain \(\Omega \subset \R^d\) with boundary \(\partial\Omega\):
\begin{subequations}\label{eq:interface}
    \begin{align}
       \cL(u) = f, ~&\quad \bx \in \cup_{k=1}^K \Omega_k, \label{eq:interface1}\\[2pt]
       \jump{u} = h_{i}, ~&\quad \bx \in \Gamma_i,~i=1,2,\cdots,m, \label{eq:interface2}\\[2pt]
       \jump{\cD(u)\cdot \bn_i} = q_{i}, ~&\quad \bx \in \Gamma_i, ~ i =1,2,\cdots,m,\label{eq:interface3}\\[2pt]
       \cB u = g, ~&\quad \bx \in \partial\Omega, \label{eq:interface4}
    \end{align}
\end{subequations}
where the domain \(\Omega\) is partitioned into \(K\) open subdomains \(\{\Omega_1,\Omega_2,\dots,\Omega_K\}\) by a collection of closed interfaces $\{\Gamma_i,~i=1,2,\cdots, m\mid \Gamma_i = \overline{\Omega}_k\cap\overline{\Omega}_l,1\le k,l\le K~\text{and}~k \neq l\}$. 
The operators $\cL$ and $\cD$ denote elliptic differential operators defined in the subdomains and on the interfaces, respectively. 
The vector $\bn_i$ denotes the unit outward normal vector on \( \Gamma_i \) from \( \Omega_k \) to \( \Omega_l \). 
For a quantity $v$, the jump across the interface $\Gamma_i$ is defined as
\begin{equation}
\jump{v} = v|_{\Omega_k} - v|_{\Omega_l}.
\end{equation}
The functions $f$, $g$, $h_i$, and $q_i$ represent the source term, the boundary data, and the two types of interface jump conditions, respectively. As an illustrative example, the bounded domain $\Omega$ shown in Figure~\ref{fig:interface_sketch} (a) is divided into $K=3$ subdomains $\{\Omega_i\}_{i=1}^3$, with interfaces $\Gamma_1$ and $\Gamma_2$.

Next we propose a novel LD-GTransNet within the lifting framework. The key idea is to lift a piecewise smooth function $u(\bx)$ into a high-dimensional space via a
newly introduced auxiliary variable $\xi$, while decoupling the spatial coordinates $\bx$ from this auxiliary variable $\xi$ so that they are processed in separate subnetworks.
{The whole architecture consists of two parallel subnetworks, namely {\em the spatial subnetwork} and {\em the auxiliary subnetwork}. The former is a GTransNet that processes spatial coordinates $\bx$ and constitutes the spatial basis functions $\{\psi_i(\bx)\}_{i=1}^N$, while the latter is a TransNet that processes the auxiliary variable $\xi$ and constitutes  the auxiliary basis functions $\{\varphi_i(\xi)\}_{i=1}^N$. The two sets of features are combined via a Hadamard product, enabling the auxiliary features to selectively gate the spatial features according to the subdomain label. The essence of this design lies in decoupling -- the spatial subnetwork learns a high-level smooth solution structures over the entire domain, while the auxiliary subnetwork acts as a region-aware switch that remains fixed during spatial differentiation.}

Specifically, the auxiliary variable $\xi=\xi(\bx)$ serves as a piecewise-constant regional indicator, assigning distinct constants $\{c_k\}$ to open subdomains $\{\Omega_k\}$. Since $\xi$ is constant on each $\Omega_k$, we have $\nabla\xi=\mathbf{0}$ everywhere on $\cup_{k=1}^K \Omega_k$; the interface $\Gamma$, being a set of measure zero, is treated separately through the jump conditions \eqref{eq:interface2}--\eqref{eq:interface3}. Consequently, all differential operators in \eqref{eq:interface} act solely on the spatial basis functions $\{\psi_i(\bx)\}_{i=1}^N$, whereas the auxiliary basis functions $\{\varphi_i(\xi)\}_{i=1}^N$ behaves as a fixed coefficient with respect to spatial variables. Interface discontinuities are not learned implicitly through the function approximation capability of a neural network, but are explicitly constructed through the difference $\varphi_i(c_k)-\varphi_i(c_l)$ of auxiliary basis functions.

\begin{figure}[ht!]
    \centering
    \includegraphics[width=1\linewidth]{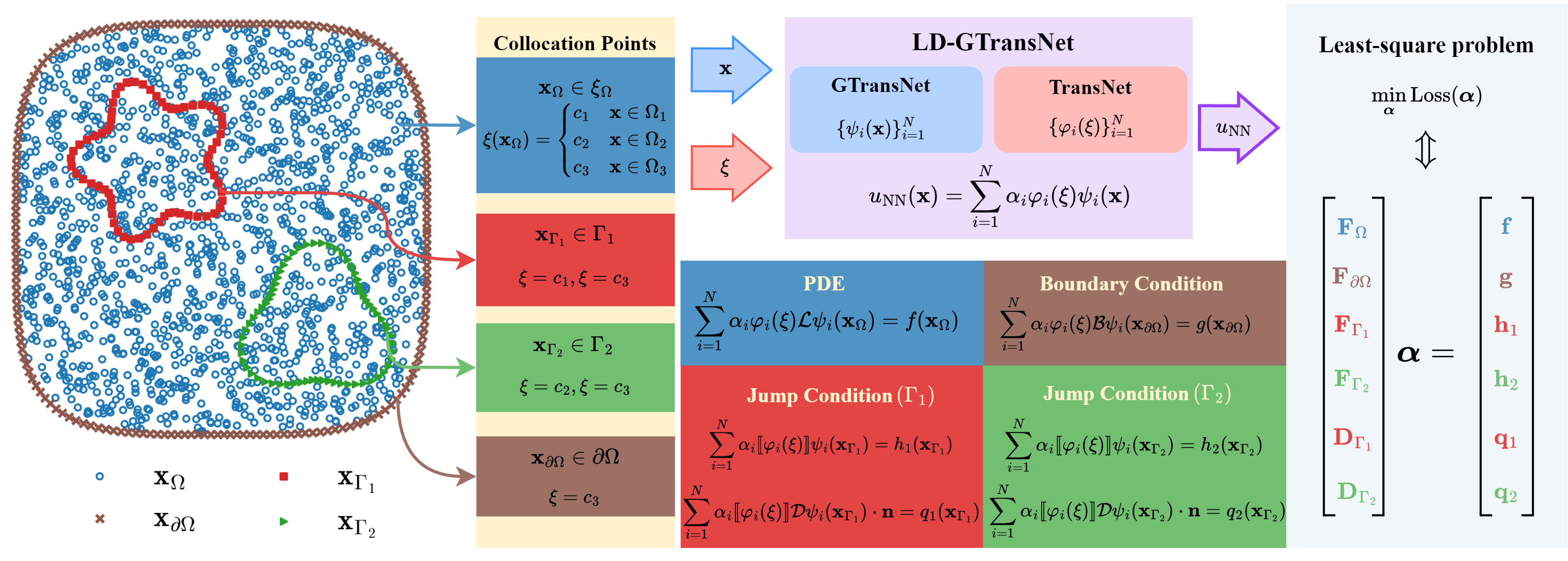}
    \caption{The {pipeline} of the proposed LD-GTransNet for solving elliptic interface {problems}.}
    \label{fig:AGTransNet}
\end{figure}

Taking the configuration shown in Figure \ref{fig:interface_sketch} as an illustrative example, the neural network solution of the elliptic interface problem \eqref{eq:interface} is  given by
\begin{equation}
u_{\rm NN}(\bx,\xi(\bx)) = \sum_{i=1}^N \alpha_i\varphi_i(\xi(\bx)) \psi_i(\bx) := \valpha^\top (\vvarphi(\xi(\bx))\odot \vpsi(\bx)),
\end{equation}
where $\valpha=(\alpha_1,\cdots,\alpha_N)^\top$, $\vvarphi(\xi)=(\varphi_1(\xi),\cdots,\varphi_N(\xi))^\top$, $\vpsi(\bx)=(\psi_1(\bx),\cdots,\psi_N(\bx))^\top$, $\xi(\bx)$ is a piecewise constant function defined by:
\begin{equation}
\xi(\bx) = \begin{cases}
    c_1, & \bx \in \Omega_1, \\
    c_2, & \bx \in \Omega_2, \\
    c_3, & \bx \in \Omega_3,
\end{cases}
\end{equation}
with $\{c_k\}$ being user-defined constants, and $\odot$ denotes element-wise multiplication (i.e., the Hadamard product). To predetermine the basis functions, the center $\bx_c$ and radius $R$ in GTransNet are selected to cover the domain $\Omega$, while the center $\xi_c$ and radius $R_\xi$ in TransNet are chosen such that their support spans the interval $[\min\{c_k\}, \max\{c_k\}]$. The corresponding pipeline of the  LD-GTransNet for solving elliptic interface problems is illustrated in Figure \ref{fig:AGTransNet}.

Assume the sets of collocation points on the domain $\Omega$, the boundary $\partial\Omega$, and the interfaces $\Gamma_1$ and $\Gamma_2$ are denoted by $\{\bx_{\Omega}^i\}_{i=1}^{K_{\Omega}}$, $\{\bx_{\partial\Omega}^i\}_{i=1}^{K_{\partial \Omega}}$, $\{\bx_{\Gamma_1}^i\}_{i=1}^{K_{\Gamma_1}}$, and $\{\bx_{\Gamma_2}^i\}_{i=1}^{K_{\Gamma_2}}$, respectively. {For brevity, we denote $\xi_{\Omega}^i = \xi(\bx_{\Omega}^i)$, which identifies the subdomain containing the collocation point $\bx_{\Omega}^i$.} The LD-GTransNet method for solving the elliptic interface problem \eqref{eq:interface} is formulated as the following minimization problem:
\begin{equation}\label{eq:optim}
\min_{\valpha}\, \text{Loss}(\valpha)
\end{equation}
where the loss function is defined as the weighted sum of the PDE residuals
\begin{equation}\label{eq:optim_equivalent}
    \begin{aligned}
        \text{Loss}( \valpha) = 
        & \sum_{i=1}^{K_{\Omega}} \lambda_1\left(\mathcal{L}(u_{\rm NN})(\bx_{\Omega}^i,\xi_{\Omega}^i) - f(\bx_{\Omega}^i)\right)^2
        + \sum_{i=1}^{K_{\partial \Omega}} \lambda_2\left(\mathcal{B}(u_{\rm NN})(\bx_{\partial\Omega}^i,c_3) - g(\bx_{\partial\Omega}^i)\right)^2 \\
        \quad & + \sum_{i=1}^{K_{\Gamma_1}} \lambda_3\left(\jump{u_{\rm NN}(\bx_{\Gamma_1}^i,\bullet)} - h_1(\bx_{\Gamma_1}^i)\right)^2 + \sum_{i=1}^{K_{\Gamma_1}} \lambda_4\left(\jump{\cD(u_{\rm NN})(\bx_{\Gamma_1}^i,\bullet)}\cdot \bn_1 - q_1(\bx_{\Gamma_1}^i)\right)^2 \\
        \quad & + \sum_{i=1}^{K_{\Gamma_2}} \lambda_5\left(\jump{u_{\rm NN}(\bx_{\Gamma_2}^i,\bullet)} - h_2(\bx_{\Gamma_2}^i)\right)^2
        + \sum_{i=1}^{K_{\Gamma_2}} \lambda_6\left(\jump{\cD(u_{\rm NN})(\bx_{\Gamma_2}^i,\bullet)}\cdot \bn_2 - q_2(\bx_{\Gamma_2}^i)\right)^2
    \end{aligned}
\end{equation}
with $\{\lambda_i\}_{i=1}^6$ being some positive balancing parameters.
Here, the jump terms are defined as 
$$
\begin{array}{rcl}
\jump{u_{\rm NN}(\bx_{\Gamma_1}^i,\bullet)} &=& u_{\rm NN}(\bx_{\Gamma_1}^i,c_1)-u_{\rm NN}(\bx_{\Gamma_1}^i,c_3),  \\
\jump{u_{\rm NN}(\bx_{\Gamma_2}^i,\bullet)} &=& u_{\rm NN}(\bx_{\Gamma_2}^i,c_2)-u_{\rm NN}(\bx_{\Gamma_2}^i,c_3), \\
\jump{\cD(u_{\rm NN})(\bx_{\Gamma_1}^i,\bullet)} &=& \cD(u_{\rm NN})(\bx_{\Gamma_1}^i,c_1)-\cD(u_{\rm NN})(\bx_{\Gamma_1}^i,c_3),\\
\jump{\cD(u_{\rm NN})(\bx_{\Gamma_2}^i,\bullet)} &=& \cD(u_{\rm NN})(\bx_{\Gamma_2}^i,c_2)-\cD(u_{\rm NN})(\bx_{\Gamma_2}^i,c_3).
\end{array}
$$

Solving the above optimization problem \eqref{eq:optim} is equivalent to solving the following linear least-squares system:
\begin{equation}\label{eq:lsq}
    \begin{bmatrix}
        \sqrt{\lambda_1} \mathbf{F}_{\Omega} \\
        \sqrt{\lambda_2} \mathbf{F}_{\partial \Omega}\\
        \sqrt{\lambda_3} \mathbf{F}_{\Gamma_1} \\        
        \sqrt{\lambda_5} \mathbf{F}_{\Gamma_2} \\
        \sqrt{\lambda_4} \mathbf{D}_{\Gamma_1} \\
        \sqrt{\lambda_6} \mathbf{D}_{\Gamma_2} \\
    \end{bmatrix}\valpha = 
    \begin{bmatrix}
        \sqrt{\lambda_1} \mathbf{F} \\
        \sqrt{\lambda_2} \mathbf{g} \\
        \sqrt{\lambda_3} \mathbf{h}_1 \\
        \sqrt{\lambda_5} \mathbf{h}_2 \\
        \sqrt{\lambda_4} \mathbf{q}_1 \\        
        \sqrt{\lambda_6} \mathbf{q}_2 \\
    \end{bmatrix},
\end{equation}
where the block matrices are formed by the Hadamard product between auxiliary-feature matrices and spatial-feature matrices:
\begin{equation}
\begin{aligned}
&\mathbf{F}_\Omega = \vvarphi_\Omega \odot \cL(\vpsi)_\Omega,\quad
\mathbf{F}_{\partial \Omega} = \vvarphi_{\partial\Omega} \odot \cB(\vpsi)_{\partial \Omega}, \\
&\mathbf{F}_{\Gamma_1} = \delta\vvarphi_{\Gamma_1} \odot \vpsi_{\Gamma_1}, \quad
\mathbf{F}_{\Gamma_2} = \delta\vvarphi_{\Gamma_2} \odot \vpsi_{\Gamma_2},\\
&\mathbf{D}_{\Gamma_1} = \delta\vvarphi_{\Gamma_1} \odot \cD(\vpsi)_{\Gamma_1}, \quad
\mathbf{D}_{\Gamma_2} = \delta\vvarphi_{\Gamma_2} \odot \cD(\vpsi)_{\Gamma_2}.
\end{aligned}
\end{equation}

The auxiliary-feature matrices are constructed by replicating the respective row vectors, which can be compactly expressed using the Kronecker product $\otimes$  with an all-ones column vector $\mathbf{1}_K=(1,\cdots,1)^\top\in\mathbb R^K$: 
\begin{equation*}
\begin{aligned}
&\vvarphi_{\Omega} = \mathbf{1}_{K_{\Omega}} \otimes 
    [\varphi_1(\xi_{\Omega}^i), \cdots, \varphi_N(\xi_{\Omega}^i) ],\\
&\delta\vvarphi_{\Gamma_1} = \mathbf{1}_{K_{\Gamma_1}} \otimes 
    [\varphi_1(c_1)-\varphi_1(c_3),\; \cdots,\; \varphi_N(c_1)-\varphi_N(c_3)],\\
&\vvarphi_{\partial\Omega} = \mathbf{1}_{K_{\partial\Omega}} \otimes 
    [\varphi_1(c_3), \cdots, \varphi_N(c_3) ]\\
&\delta\vvarphi_{\Gamma_2} = \mathbf{1}_{K_{\Gamma_2}} \otimes 
    [\varphi_1(c_2)-\varphi_1(c_3),\; \cdots,\; \varphi_N(c_2)-\varphi_N(c_3)].
\end{aligned}
\end{equation*}
Here  $\delta \vvarphi_{\Gamma_i}$ 
signifies the auxiliary-feature difference across the interface $\Gamma_i$. 
The spatial-feature matrices are given by
$$
\cL(\vpsi)_\Omega = \begin{bmatrix}
        \cL(\psi_1)(\bx_{\Omega}^1) &  
        \cdots & \cL(\psi_N)(\bx_{\Omega}^1) \\ 
        \vdots & \ddots & \vdots \\
       \cL(\psi_1)(\bx_{\Omega}^{K_\Omega})   &
        \cdots & \cL(\psi_N)(\bx_{\Omega}^{K_\Omega}) \\
    \end{bmatrix},
$$
$$
\cB(\vpsi)_{\partial \Omega} = \begin{bmatrix}
        \cB(\psi_1)(\bx_{\partial \Omega}^1)  &
        \cdots & \cB(\psi_N)(\bx_{\partial \Omega}^1) \\
        \vdots & \ddots & \vdots \\
       \cB(\psi_1)(\bx_{\partial \Omega}^{K_{\partial \Omega}})  &
        \cdots & \cB(\psi_N)(\bx_{\partial \Omega}^{K_{\partial \Omega}}) \\
    \end{bmatrix},
$$

$$
\vpsi_{\Gamma_1} = \begin{bmatrix}
        \psi_1(\bx_{\Gamma_1}^1)  &
        \cdots & \psi_N(\bx_{\Gamma_1}^1) \\
        \vdots & \ddots & \vdots \\
       \psi_1(\bx_{\Gamma_1}^{K_{\Gamma_1}})  &
        \cdots & \psi_N(\bx_{\Gamma_1}^{K_{\Gamma_1}}) 
    \end{bmatrix}, \quad
\vpsi_{\Gamma_2} = \begin{bmatrix}
        \psi_1(\bx_{\Gamma_2}^1)  &
        \cdots & \psi_N(\bx_{\Gamma_2}^1) \\
        \vdots & \ddots & \vdots \\
       \psi_1(\bx_{\Gamma_2}^{K_{\Gamma_2}})  &
        \cdots & \psi_N(\bx_{\Gamma_2}^{K_{\Gamma_2}}) 
    \end{bmatrix}, 
    $$
$$
\cD(\vpsi)_{\Gamma_1}=\begin{bmatrix}
        \cD(\psi_1)(\bx_{\Gamma_1}^1) \cdot \bn_1 & \cdots & \cD(\psi_N)(\bx_{\Gamma_1}^1) \cdot \bn_1 \\
        \vdots & \ddots & \vdots \\
        \cD(\psi_1)(\bx_{\Gamma_1}^{K_{\Gamma_1}}) \cdot \bn_1 & \cdots & \cD(\psi_N)(\bx_{\Gamma_1}^{K_{\Gamma_1}}) \cdot \bn_1 \\
    \end{bmatrix}, 
$$
$$
\cD(\vpsi)_{\Gamma_2}=\begin{bmatrix}
        \cD(\psi_1)(\bx_{\Gamma_2}^1) \cdot \bn_2 & \cdots & \cD(\psi_N)(\bx_{\Gamma_2}^1) \cdot \bn_2 \\
        \vdots & \ddots & \vdots \\
        \cD(\psi_1)(\bx_{\Gamma_2}^{K_{\Gamma_2}}) \cdot \bn_2 & \cdots & \cD(\psi_N)(\bx_{\Gamma_2}^{K_{\Gamma_2}}) \cdot \bn_2 \\
    \end{bmatrix}.
$$
Finally, the right-hand side vectors are defined by
$$
    \mathbf{F} = \begin{bmatrix}
        f(\bx_{\Omega}^1) \\
        \vdots \\
        f(\bx_{\Omega}^{K_\Omega}) 
    \end{bmatrix},\quad
    \mathbf{g} = \begin{bmatrix}
        g(\bx_{\partial\Omega}^1) \\
        \vdots \\
        g(\bx_{\partial\Omega}^{K_{\partial \Omega}}) 
    \end{bmatrix}, \quad
    \mathbf{h}_i = \begin{bmatrix}
        h_i(\bx_{\Gamma_i}^1) \\
        \vdots \\
        h_i(\bx_{\Gamma_i}^{K_{\Gamma_i}}) 
    \end{bmatrix},\quad
    \mathbf{q}_i = \begin{bmatrix}
        q_i(\bx_{\Gamma_i}^1) \\
        \vdots \\
        q_i(\bx_{\Gamma_i}^{K_{\Gamma_i}}) 
    \end{bmatrix},\quad
    i=1,2.
$$

{The key consequence of this separation is that the PDE residuals in the subdomain interiors and the interface jump residuals are decoupled into mutually independent block rows in the least-squares system \eqref{eq:lsq}.} Specifically, the PDE blocks involve the Hadamard product $\varphi_\Omega \odot \mathcal{L}(\psi)_\Omega$, while the interface blocks directly encode jumps through the difference $\varphi(c_k)-\varphi(c_l)$ of auxiliary basis functions. This structure circumvents the inherent difficulty of approximating discontinuous functions with a single neural network, while replacing the nonlinear optimization in conventional methods with a linear least-squares formulation. All residual terms are incorporated into a unified least-squares framework, with the penalty parameters automatically balancing the magnitudes of the corresponding block rows through row scaling:
\begin{equation}
    \begin{aligned}
        \sqrt\lambda_1 &= \frac{1}{\max_{1\le i \le K_\Omega} \left(\max_{1\le j \le N} |(\mathbf{F}_\Omega)_{ij}|, |f(\bx_{\Omega}^i)|\right)}, \\
        \sqrt\lambda_2 &= \frac{1}{\max_{1\le i \le K_{\partial \Omega}} \left(\max_{1\le j \le N} |(\mathbf{F}_{\partial \Omega})_{ij}|, |g(\bx_{\partial\Omega}^i)|\right)}, \\
        \sqrt\lambda_3 &= \frac{1}{\max_{1\le i \le K_{\Gamma_1}} \left(\max_{1\le j \le N} |(\mathbf{F}_{\Gamma_1})_{ij}|, |h_1(\bx_{\Gamma_1}^i)|\right)}, \\
        \sqrt\lambda_4 &= \frac{1}{\max_{1\le i \le K_{\Gamma_1}} \left(\max_{1\le j \le N} |(\mathbf{D}_{\Gamma_1})_{ij}|, |q_1(\bx_{\Gamma_1}^i)|\right)}.
    \end{aligned}
\end{equation}
 {The above selection of penalty parameters balances the magnitudes of the residual terms in the loss function \eqref{eq:optim_equivalent}, preventing any single residual from dominating the objective function during optimization. This is an empirically motivated heuristic strategy \cite{Lu2025} that has been observed to improve the conditioning of the least-squares matrix for interface problems; however, a rigorous theoretical analysis of its convergence properties remains to be investigated.} The selection of $\lambda_5, \lambda_6$ follows the same principle as $\lambda_3, \lambda_4$ and is omitted here.

{Moreover, LD-GTransNet remains a natural and strict generalization of GTransNet.} In the absence of interfaces, all points in the domain $\Omega$ share the same subdomain label, say $c_1$. Consequently, for any pair of collocation points on a nominal interface, the auxiliary-variable difference terms vanish because $\varphi_i(c_1) - \varphi_i(c_1) = 0$ for all $i$. Hence $\delta\varphi_{\Gamma_1}$ and $\delta\varphi_{\Gamma_2}$ become zero matrices, which in turn implies that $\mathbf{F}_{\Gamma_1}$, $\mathbf{F}_{\Gamma_2}$, $\mathbf{D}_{\Gamma_1}$, and $\mathbf{D}_{\Gamma_2}$ are also zero matrices. The problem ~\eqref{eq:lsq} then reduces to a system involving only its first two block rows, which corresponds exactly to the least-squares formulation of the original GTransNet method.

\subsection{Moving interface problems}\label{sec:3.2}
In this subsection, we extend the LD-GTransNet method investigated in the preceding subsection to time-dependent moving interface problems. By treating time $t$ as an additional coordinate, the moving interface $\Gamma(t)$ in physical space becomes a stationary manifold in the space-time domain $\Omega_T := \Omega \times (0,T]$, while the physical boundary remains fixed, yielding the space-time cylinder $\partial\Omega\times(0,T]$.
If the interface motion is known a priori, it suffices to include $t$ as an input to GTransNet, and the problem reduces to the static case described in Section~3.1. When the interface motion is unknown and coupled with the solution, the problem becomes nonlinear; we decouple it via Gauss-Seidel block iteration (i.e., block coordinate descent).

In the following we illustrate our method for a generic moving interface problem with two phases; its extension to multiple interfaces and phases is straightforward. Assume that the moving interface is defined as the zero level-set of an unknown function $\zeta(\bx,t)$:
\begin{equation}\label{eq:interface_def}
    \Gamma(t) = \{\bx \mid \zeta(\bx,t) = 0\},
\end{equation}
which partitions $\Omega$ into two evolving subdomains $\Omega_1(t)$ and $\Omega_2(t)$. The interface motion is governed by a solution-dependent condition of the general form
\begin{equation}\label{eq:interface_condition}
    \mathcal{F}(\zeta, u) = 0, \quad (\bx, t) \in \Gamma(t),
\end{equation}
where $\mathcal{F}$ is a (possibly nonlinear) operator that may involve temporal and spatial derivatives of $\zeta$ as well as derivatives of $u$. The governing PDE in each subdomain, together with the interface jump conditions and initial/boundary data, reads:
\begin{subequations}\label{eq:general_moving}
    \begin{align}
        \dfrac{\partial u}{\partial t} - \cL(u) = f, & \quad (\bx,t) \in \Omega_1(t) \cup \Omega_2(t), \\[6pt]
        \jump{u} = h, &\quad  (\bx,t) \in \Gamma(t), \\[4pt]
        \jump{\cD(u)\cdot \bn} = q, &\quad  (\bx,t) \in \Gamma(t), \\[4pt]
        \cB u = g, &\quad  (\bx,t) \in \partial\Omega\times(0,T], \\[4pt]
        u(\bx,0) = u_0, &\quad  \bx \in \Omega_1(0) \cup \Omega_2(0).
    \end{align}
\end{subequations}
We remark that the equations~\eqref{eq:interface_condition} and~\eqref{eq:general_moving} form a coupled nonlinear system: $u$ depends on the interface location $\Gamma(t)$ (determined by $\zeta$), while $\zeta$ is constrained by $\mathcal{F}(\zeta,u)=0$.

To solve this coupled system, we approximate $u$ with the LD-GTransNet method as follows:
\begin{equation}
u_{\rm NN}(\bx,t,\xi) = \sum_{i=1}^{N_u} \alpha_i\; \varphi_i(\xi(\bx,t))\,\psi_i(\bx,t) = \valpha^\top \bigl(\vvarphi(\xi(\bx,t)) \odot \vpsi(\bx,t)\bigr),
\end{equation}
and approximate $\zeta$ by a separate TransNet through
\begin{equation}\label{znn}
\zeta_{\rm NN}(\bx,t) = \sum_{i=1}^{N_\zeta} \beta_i \phi_i(\bx,t) = \boldsymbol{\beta}^\top \boldsymbol{\phi}(\bx,t).
\end{equation}
Here the auxiliary variable $\xi$ is defined as a piecewise-constant function of the level-set value: $\xi = c_1$ if $\zeta(\bx,t) < 0$ and $\xi = c_2$ if $\zeta(\bx,t) > 0$, with the interface $\Gamma(t)$ corresponding to the jump discontinuity.
The corresponding optimization problem is then formulated as minimizing the following joint loss function:
\begin{equation}\label{eq:loss_joint}
    \min_{(\valpha, \boldsymbol{\beta})} \mathrm{Loss}(\valpha, \boldsymbol{\beta}) = \mathrm{Loss}_u(\valpha, \boldsymbol{\beta}) + \mathrm{Loss}_\zeta(\valpha, \boldsymbol{\beta}),
\end{equation}
where $\mathrm{Loss}_u$ follows the same construction as in the elliptic case~\eqref{eq:optim_equivalent}, but with collocation points now distributed in the space-time domain $\Omega_T$ and supplemented by the initial condition term:
\begin{equation}\label{eq:loss_u_moving}
    \begin{aligned}
        \mathrm{Loss}_u(\valpha,\boldsymbol{\beta}) = 
        & \sum_{i=1}^{K_{\Omega_T}} \lambda_1\left(\frac{\partial u_{\rm NN}}{\partial t} - \mathcal{L}(u_{\rm NN})(\bx_{\Omega}^i,t^i,\xi_{\Omega_T}^i) - f(\bx_{\Omega}^i,t^i)\right)^2 \\
        & + \sum_{i=1}^{K_{\partial\Omega_T}} \lambda_2\left(\mathcal{B}(u_{\rm NN})(\bx_{\partial\Omega}^i,t^i,c_2) - g(\bx_{\partial\Omega}^i,t^i)\right)^2 \\
        & + \sum_{i=1}^{K_{\Gamma}} \lambda_3\left(\jump{u_{\rm NN}(\bx_{\Gamma}^i,\bullet)} - h(\bx_{\Gamma}^i,t^i)\right)^2 \\
        & + \sum_{i=1}^{K_{\Gamma}} \lambda_4\left(\jump{\cD(u_{\rm NN})(\bx_{\Gamma}^i,\bullet)}\cdot \bn - q(\bx_{\Gamma}^i,t^i)\right)^2 \\
        & + \sum_{i=1}^{K_{\Omega_0}} \lambda_0\left(u_{\rm NN}(\bx_{\Omega_0}^i,0,\xi_{\Omega_0}^i) - u_0(\bx_{\Omega_0}^i)\right)^2,
    \end{aligned}
\end{equation}
where $\Gamma$ denotes the interface determined by the current $\zeta_{\rm NN}$, $\Omega_0 = \Omega_1(0) \cup \Omega_2(0)$ is the domain at initial time, $\xi_{\Omega_T}^i = \xi(\bx_{\Omega}^i,t^i)$ identifies the phase of the space-time collocation point, and $\{\lambda_i\}_{i=0}^{4}$ are positive balancing parameters.
The interface loss $\mathrm{Loss}_\zeta$ is defined as:
\begin{equation}\label{eq:loss_zeta}
    \mathrm{Loss}_\zeta(\valpha,\boldsymbol{\beta}) =
    \sum_{i=1}^{K_{\Gamma}} \mu_1\left(\mathcal{F}(\zeta_{\rm NN}, u_{\rm NN})(\bx_{\Gamma}^i,t^i)\right)^2
    + \sum_{i=1}^{K_{\Omega_0}} \mu_2\left(\zeta_{\rm NN}(\bx_{\Omega_0}^i,0) - \zeta_0(\bx_{\Omega_0}^i)\right)^2,
\end{equation}
where $\mu_1$ and $\mu_2$ are again two positive balancing parameters.  

Since the loss function~\eqref{eq:loss_joint} is jointly nonconvex in $(\valpha, \boldsymbol{\beta})$, direct optimization is challenging. We therefore employ Gauss-Seidel block iteration to decouple the problem into two alternating linear least-squares subproblems: solve
\begin{equation}\label{eq:gauss_seidel}
    \begin{cases}
        \valpha^{n+1} = \arg\min_{\valpha} \mathrm{Loss}_u(\valpha,\boldsymbol{\beta}^n), \\[6pt]
        \boldsymbol{\beta}^{n+1} = \arg\min_{\boldsymbol{\beta}} \mathrm{Loss}_\zeta(\valpha^{n+1},\boldsymbol{\beta}),
    \end{cases}
\quad n = 0,1,2,\dots,
\end{equation}
until the approximate solutions converge.
At the $n$-th iteration, we first solve the linear least-squares problem for $\valpha$ with fixed $\boldsymbol{\beta}^n$, obtaining $\valpha^{n+1}$; we then immediately use the updated $\valpha^{n+1}$ to solve the linear least-squares problem for $\boldsymbol{\beta}$. Since both subnetworks employ fixed hidden-layer parameters, each subproblem minimization reduces to a linear least-squares problem: $\mathrm{Loss}_u$ is quadratic in $\valpha$, and $\mathrm{Loss}_\zeta$ is quadratic in $\boldsymbol{\beta}$.

\begin{remark}
The classical two-phase Stefan problem in two dimensions serves as a concrete example of the above framework. For an interface represented by $x=s(y,t)$, the two subdomains are $\Omega_1(t)=\{(x,y)\in\Omega\mid x<s(y,t)\}$ and $\Omega_2(t)=\{(x,y)\in\Omega\mid x>s(y,t)\}$. The diffusion coefficient is $\kappa=\kappa_1$ in $\Omega_1(t)$ and $\kappa=\kappa_2$ in $\Omega_2(t)$, where $\kappa_1,\kappa_2>0$ are constants. The governing PDE is given by
\begin{equation}\label{eq:stefan}
\dfrac{\partial u}{\partial t} - \kappa \Delta u = 0, \quad (x,y,t) \in \Omega_1(t) \cup \Omega_2(t),
\end{equation}
with jump conditions $\jump{u} = h,~\jump{\kappa\nabla u\cdot \bn} = q$ on $\Gamma(t)$ and appropriate initial/boundary data. The interface is parameterized as $\Gamma(t) = \{(x,y,t)\mid x = s(y,t)\}$, which corresponds to choosing the level-set function $\zeta(x,y,t) = x - s(y,t)$ in~\eqref{eq:interface_def}.
The operator $\mathcal{F}$ in~\eqref{eq:interface_condition} then takes the specific form of the Stefan condition:
\begin{equation}
\mathcal{F}(\zeta,u) = \partial_t \zeta + \jump{\kappa\nabla u\cdot\bn} = 0,
\end{equation}
 where $\jump{\kappa\nabla u\cdot\bn}$ denotes the jump of normal heat flux across the interface. 
Alternatively, we also can directly represent $s(y,t)$ by a TransNet (instead of $\zeta_{\rm NN}$ for $\zeta(x,y,t)$)
\begin{equation}
s_{\rm NN}(y,t) = \sum_{i=1}^{N_s} \beta_i \phi_i(y,t) = \boldsymbol{\beta}^\top \boldsymbol{\phi}(y,t),
\end{equation}
then the interface loss~\eqref{eq:loss_zeta} consequently reduces to
  \begin{equation}
\mathrm{Loss}_s(\valpha,\boldsymbol{\beta})  = \sum_{i=1}^{K_{\Gamma^n}} \mu_1\left(\partial_t s_{\rm NN}(y^i,t^i) + \jump{\kappa\nabla u_{\rm NN}\cdot\bn}(y^i,t^i)\right)^2 + \sum_{i=1}^{K_{\Omega_0}} \mu_2\left(s_{\rm NN}(y^i,0) - s_0(y^i)\right)^2.
 \end{equation}
\end{remark}
\section{Numerical Experiments}\label{sec:4}
In this section, we conduct a series of numerical experiments in both two and three dimensions to validate the accuracy, robustness, and efficiency of the LD-GTransNet method. A variety of interface problems are considered, including a two-dimensional isotropic diffusion equation with multiple interfaces, a two-dimensional Stokes problem with a circular interface, a two-dimensional anisotropic diffusion-reaction equation involving multiple interfaces, a three-dimensional elasticity problem with an elliptic interface and two moving interface problems. 
In all experiments, the number of neurons in the last hidden layer of our LD-GTransNet is denoted by $N$. Unless stated otherwise, the shape parameter is set to $\gamma = 1$. Training points are generated using the Latin hypercube sampling method \cite{Latin}.



To evaluate the numerical accuracy, we introduce the relative $L^2$ error and relative $L^\infty$ error, defined as follows:
\begin{equation}
err_{L^{\infty}} =\frac{\|u - u_{\mathrm{NN}}\|_{L^\infty}}{\|u\|_{L^\infty}},
\quad
err_{L^2} = \frac{\|u - u_{\mathrm{NN}}\|_{L^2}}{\|u\|_{L^2}} 
\end{equation}
where $u_{\mathrm{NN}}$ denotes the numerical solution obtained via the LD-GTransNet or other methods, and $u$ is the corresponding exact solution. For 2D and 3D examples,
we use $\mathcal{O}(10^4)$ and $\mathcal{O}(10^5)$ test points, respectively, to calculate the error. All experiments are performed on a Windows system equipped with an Intel Core i5-13600KF processor. SciPy's \verb|scipy.linalg.lstsq| API is used to solve the least squares problems. All results are obtained by repeating each experiment 10 times.

\subsection{Tests for elliptic interface problems}
In this subsection, we test the proposed LD-GTransNet for some elliptic interface problems and compare it with existing methods to demonstrate its superior performance. 

\subsubsection{A 2D isotropic diffusion multiple-interface equation}\label{subsec:multi-interfaces}
\begin{figure}[ht!]
    \centering
    \includegraphics[width=1\linewidth]{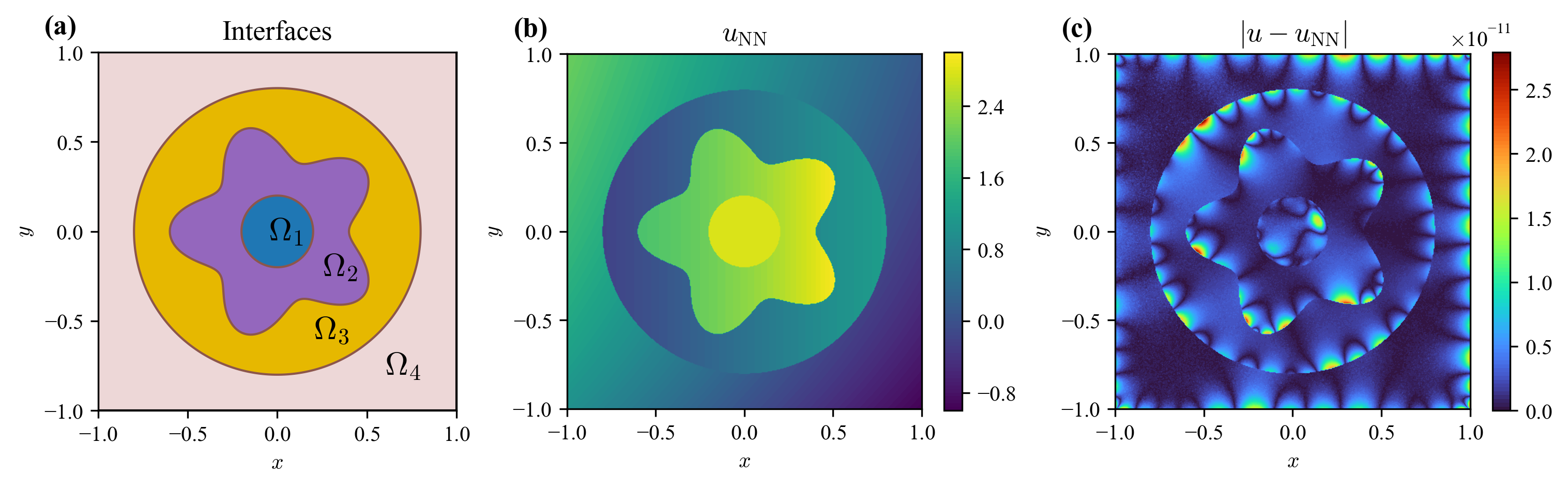}
    \vspace{-0.7cm}
    \caption{(a): Illustration of the partition of the domain into four subdomains by three interfaces. (b)-(c): Numerical solutions and pointwise absolute errors produced by the two-hidden-layer LD-GTransNet method with $N=1200$ and $(\kappa_1,\kappa_2,\kappa_3,\kappa_4) = (1,10^2,10^4,10^6)$ for the 2D diffusion {multiple-interface} problem \eqref{eq:multi-interfaces}.}
    \label{fig:multi-interface}
\end{figure}

We now consider a two-dimensional diffusion {multiple-interface}  problem on the domain $\Omega = [-1,1]^2$, featuring three interfaces as studied in \cite{Yunlong2025}. The problem is governed by the following elliptic equation:
\begin{equation}\label{eq:multi-interfaces}
\begin{cases}
-\nabla \cdot (\kappa \nabla u) = f, & (x,y) \in \bigcup_{k=1}^4\Omega_k, \\
\jump{u} = h_i, & (x,y) \in \Gamma_i,~i=1,2,3, \\
\jump{\kappa \nabla u \cdot \bn_i} = q_i, & (x,y) \in \Gamma_i,~i=1,2,3, \\
u = g, & (x,y) \in \partial \Omega,
\end{cases}
\end{equation}
where the diffusion coefficient $\kappa$ takes the constant value $\kappa_k$ in each subdomain $\Omega_k$ ($k=1,2,3,4$). As illustrated in Figure~\ref{fig:multi-interface}(a), the three interfaces are defined by
\begin{align}
\Gamma_1 &= \left\{(x,y)\mid x^2 + y^2 = 0.2^2 \right\}, \\
\Gamma_2 &= \left\{(x,y)\mid x^2 + y^2 = (0.5-0.1\cos(5\theta))^2,~\theta=\arctan\left(\frac{y}{x}\right)\in[0,2\pi)\right\}, \\
\Gamma_3 &= \left\{(x,y)\mid x^2 + y^2 = 0.8^2\right\}.
\end{align}
The exact solution to the problem \eqref{eq:multi-interfaces} is given by:
\begin{equation*}
u(x,y) =
\begin{cases}
\cos y + 1.8, & (x,y) \in \Omega_1, \\
\exp x + 1.3, & (x,y) \in \Omega_2, \\
\sin x + 0.5, & (x,y) \in \Omega_3, \\
-x +\ln(y+2), & (x,y) \in \Omega_4.
\end{cases}
\end{equation*}
The corresponding source term $f$, boundary condition $g$, and jump conditions $h_i, q_i$ ($i=1,2,3$) are determined {from the exact solution}. 
The regional indicator $\xi$ is set as 
\begin{equation*}
\xi =
\begin{cases}
1, & (x,y) \in \Omega_1, \\
0.3, & (x,y) \in \Omega_2, \\
-0.3, & (x,y) \in \Omega_3, \\
-1, & (x,y) \in \Omega_4,
\end{cases}
\end{equation*}
and the parameters are set to $\mathbf{x}_c=(0,0),R=1.5$ for GTransNet and $\xi_c=0,R_{\xi}=1.1$ for TransNet.
Collocation points are generated via Latin hypercube sampling—1500 in $\Omega$, 400 on $\partial\Omega$, and 80 on each of the interfaces $\Gamma_1,\Gamma_2,\Gamma_3$.

\begin{figure}[ht!]
    \centering
    \includegraphics[width=1\linewidth]{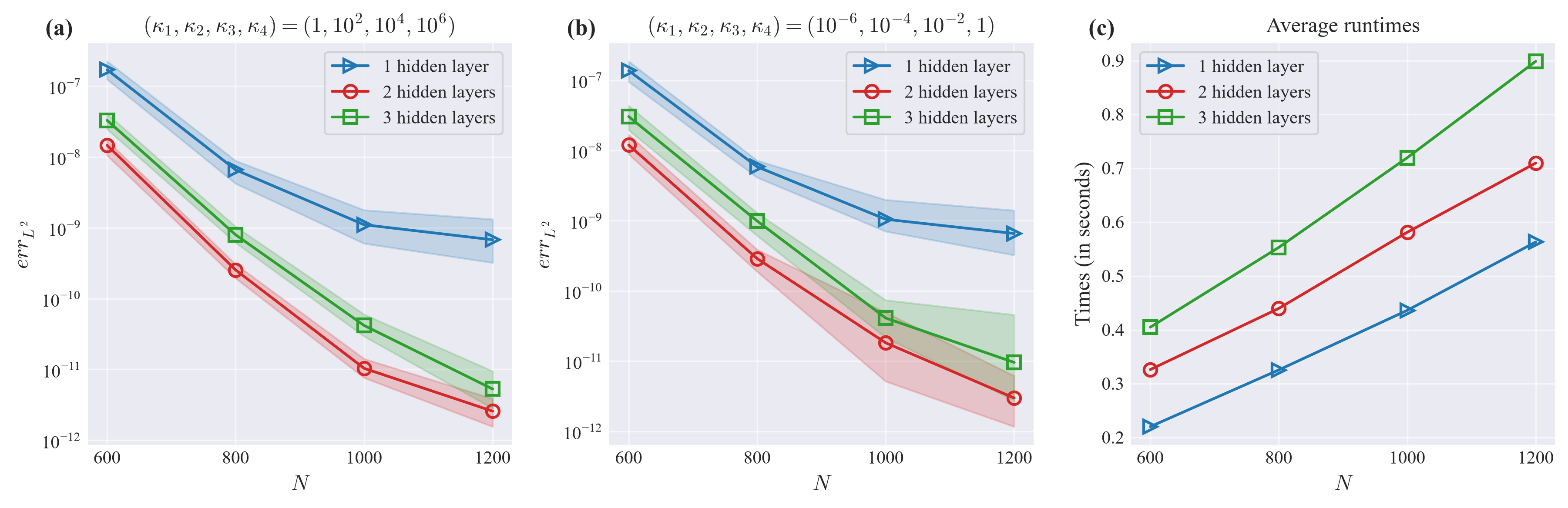}
    \vspace{-0.2cm}
    \caption{(a)-(b): Relative $L^2$ errors of the numerical solutions obtained by the LD-GTransNet method with one, two and three hidden layers with respect to the number of neurons $N$ in the last hidden layer for the 2D diffusion {multiple-interface}  problem \eqref{eq:multi-interfaces} under different combinations of diffusion coefficients. (c): Average running times (in seconds) of the LD-GTransNet method with one, two and three hidden layers.}
    \label{fig:multi-interface-errors}
\end{figure}

\begin{table}[ht!]
  \centering\small
  \caption{Performance comparison among LD-GTransNet, LRNN \cite{Yunlong2025}, and Multi-TransNet \cite{Lu2025} for the 2D diffusion {multiple-interface} problem \eqref{eq:multi-interfaces} under different combinations of diffusion coefficients. \#constraints denotes the row number of the resulting least squares system, and \#DOFs denotes the column number of the resulting least squares system.}
  \begin{tabular}{llccc}
    \toprule
    $(\kappa_1, \kappa_2, \kappa_3, \kappa_4)$ & Method & \#Constraints & \#DOFs & $err_{L^2}(u)$ \\
    \midrule
    \multirow{3}{*}{$(1, 10^2, 10^4, 10^6)$}
    & LRNN & 6500 & 1280 & 4.64e-5 \\
    & Multi-TransNet & 6500 & 800 & 2.41e-8 \\
    & LD-GTransNet & 2380 & 800 & \textbf{1.96e-10} \\
    & & 2380 & 1200 & \textbf{2.48e-12} \\
    \midrule
    \multirow{3}{*}{$(10^{-6}, 10^{-4}, 10^{-2}, 1)$}
    & LRNN & 6500 & 1280 & 4.44e-6 \\
    & Multi-TransNet & 6500 & 800 & 2.25e-8 \\
    & LD-GTransNet & 2380 & 800 & \textbf{1.21e-10} \\
    & & 2380 & 1200 & \textbf{2.43e-12} \\
    \bottomrule
  \end{tabular}
  \label{tab:comparison_multi_interface}
\end{table}

First, we investigate the effect of the number of hidden layers on the LD-GTransNet method. For $L=2$, the number of neurons in the first hidden layer is fixed at 1200; for $L=3$, the numbers of neurons in both the first and second hidden layers are fixed at 1200. Figure \ref{fig:multi-interface-errors} (a)-(b) presents the variations of the relative $L^2$ errors with respect to the number of neurons $N$ in the last layer under different numbers of hidden layers and different combinations of diffusion coefficients.

The results indicate that even with a single hidden layer ($L=1$), our LD-GTransNet method already achieves a very small error. When two or three hidden layers are used, the errors further decrease to the order of $\mathcal{O}(10^{-12})$, demonstrating that increasing the number of hidden layers can substantially enhance numerical accuracy. Moreover, the performance with two hidden layers is comparable to that with three, suggesting limited improvement beyond two layers. Figure~\ref{fig:multi-interface-errors}(c) presents the corresponding running times, showing that increasing the number of hidden layers also leads to higher computational cost.
{Considering both accuracy and computational efficiency, the two-hidden-layer configuration provides the best overall trade-off} for LD-GTransNet. Therefore, all subsequent results are based on this configuration. For illustration, Figures~\ref{fig:multi-interface}(b)–(c) display the numerical solution and pointwise absolute error obtained by the two-hidden-layer LD-GTransNet with $N=1200$ and $(\kappa_1,\kappa_2,\kappa_3,\kappa_4) = (1,10^2,10^4,10^6)$. 

Next, we compare LD-GTransNet with LRNN \cite{Yunlong2025} and Multi-TransNet  \cite{Lu2025} on the same diffusion  {multiple-interface} problem under two sets of diffusion coefficients. As shown in Table~\ref{tab:comparison_multi_interface}, LD-GTransNet not only requires a smaller-scale linear system but also achieves errors that are seven orders of magnitude lower than LRNN and four orders of magnitude lower than Multi-TransNet, demonstrating its clear superiority over these baseline methods. {The reduction in the number of constraints stems from the decoupled architecture: by processing the auxiliary variable through a separate small TransNet rather than concatenating it with spatial inputs, LD-GTransNet avoids the need for densely sampling the lifted space and instead leverages the Hadamard product to selectively activate spatial features per subdomain. This design reduces the required number of collocation points while preserving the full expressivity of the spatial basis.}

\subsubsection{A 2D Stokes interface problem}
In this subsection, we solve the following 2D Stokes interface problem \cite{Chi2024} on the domain $\Omega = [-2,2]^2$:
\begin{equation}\label{eq:Stokes}
\begin{cases}
-\mu \Delta \bu + \nabla p = \mathbf{f}, & (x,y) \in \Omega_1 \cup \Omega_2, \\
 \nabla \cdot \bu = 0, & (x,y) \in \Omega_1 \cup \Omega_2, \\
 \jump{\bu} = \boldsymbol{0}, & (x,y) \in \Gamma, \\
 \jump{\boldsymbol{\sigma}(\bu,p)\bn} = \mathbf{q}, & (x,y) \in \Gamma, \\
 \bu = \mathbf{g}, & (x,y) \in \partial \Omega,
\end{cases}
\end{equation}
where the viscosity $ \mu(x,y) = \mu_k $ for $ (x,y) \in \Omega_k $ ($k = 1,2$), the velocity vector $ \bu = (u,v)^\top $ and the stress tensor $\boldsymbol{\sigma}(\bu,p) = -p\mathbf{I} + \mu\left(\nabla \bu + \nabla \bu^\top\right) $.
The interface is defined as
\begin{equation}
\Gamma = \{(x,y)\mid x^2+y^2=1\}
\end{equation}
The subdomains are $\Omega_1=\{(x,y)\in(-2,2)^2\mid x^2+y^2<1\}$ and $\Omega_2=\{(x,y)\in(-2,2)^2\mid x^2+y^2>1\}$. The viscosity is $\mu_1$ inside the circle and $\mu_2$ outside. The exact {solution is} given by
\begin{equation*}
u(x,y) =
\begin{cases}
\dfrac{y}{4}\left(x^2 + y^2\right), & (x,y) \in \Omega_1, \\[6pt]
\dfrac{y}{\sqrt{x^2 + y^2}} - \dfrac{3y}{4}, & (x,y) \in \Omega_2,
\end{cases}
\end{equation*}
\begin{equation*}
v(x,y) =
\begin{cases}
-\dfrac{xy^2}{4}, & (x,y) \in \Omega_1, \\[6pt]
-\dfrac{x}{\sqrt{x^2 + y^2}} + \dfrac{x}{4}\left(3 + x^2\right), & (x,y) \in \Omega_2,
\end{cases}
\end{equation*}
and
\begin{equation*}
p(x,y) =
\begin{cases}
5, & (x,y) \in \Omega_1, \\[6pt]
\left(-\dfrac{3}{4}x^3 + \dfrac{3}{8}x\right)y, & (x,y) \in \Omega_2.
\end{cases}
\end{equation*}
The source term $\mathbf{f}$, boundary condition $\mathbf{g}$, and jump conditions $\mathbf{q}$ are determined {from the exact solution}. {Note that the incompressibility $\nabla\cdot\bu=0$ is enforced as an additional residual term in the least-squares formulation~\eqref{eq:optim_equivalent}, analogous to the PDE residual.}
We set \(\xi\) as
\begin{equation*}
\xi =
\begin{cases}
1, & (x,y) \in \Omega_1, \\
-1, & (x,y) \in \Omega_2,
\end{cases}
\end{equation*}
and the parameters are set to $\mathbf{x}_c=(0,0),R=2.9$ for GTransNet and $\xi_c=0,R_{\xi}=1.1$ for TransNet.
The number of neurons in the first hidden layer is fixed at 3000.
The collocation points are distributed as follows: 3000 points in $\Omega$, 800 points on $\partial\Omega$, and 200 points on the interface $\Gamma$.

\begin{figure}[ht!]
    \centering
    \includegraphics[width=0.8\linewidth]{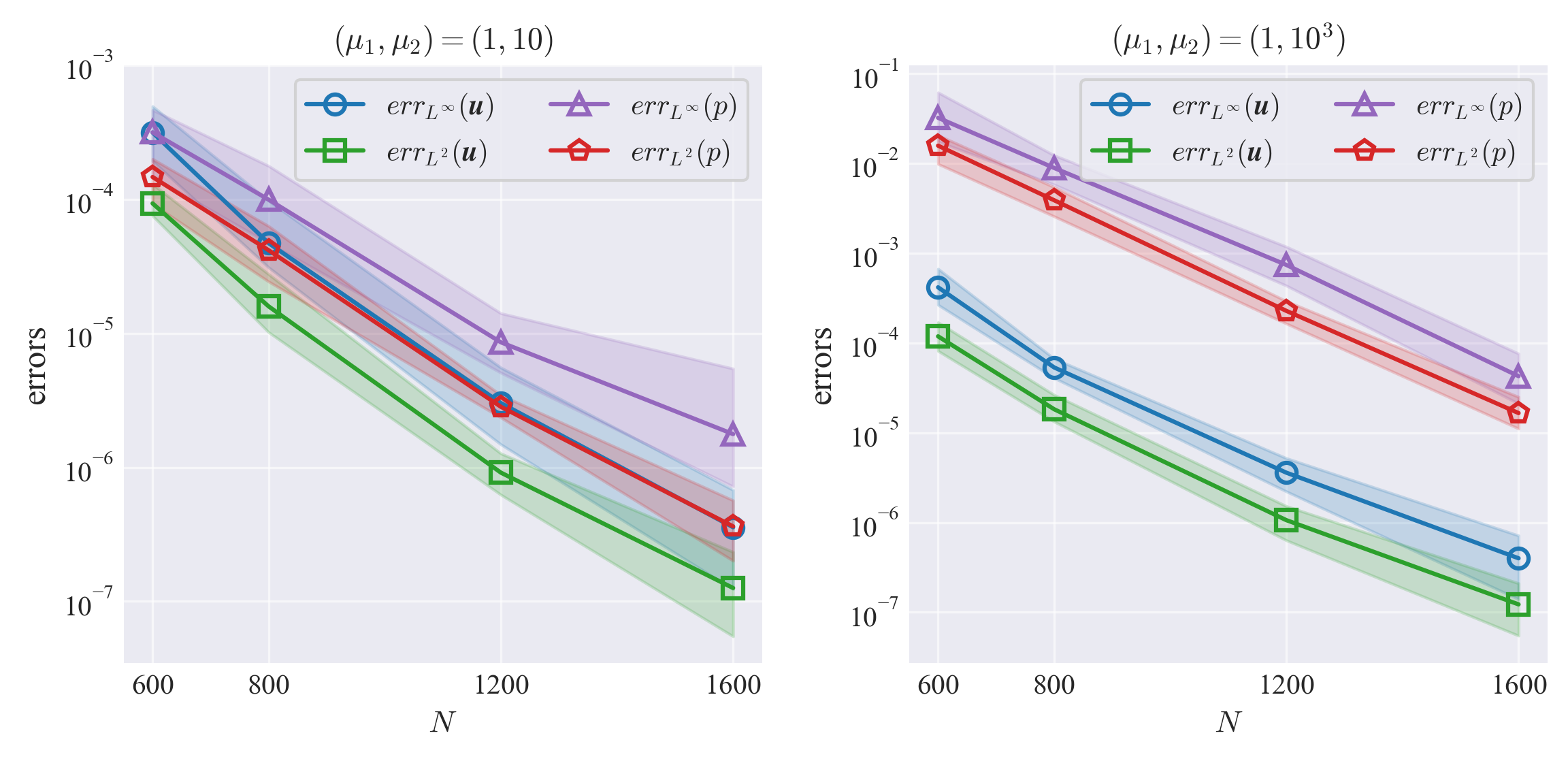}
     \vspace{-0.2cm}
    \caption{Relative $L^{\infty}$ and $L^2$ errors of the numerical solutions obtained by the LD-GTransNet method with respect to the number of neurons $N$ in the last hidden layer for the 2D Stokes interface problem \eqref{eq:Stokes} under different viscosity contrasts.}
    \label{fig:Stokes-errors}
\end{figure}

\begin{figure}[ht!]
    \centering
    \includegraphics[width=1\linewidth]{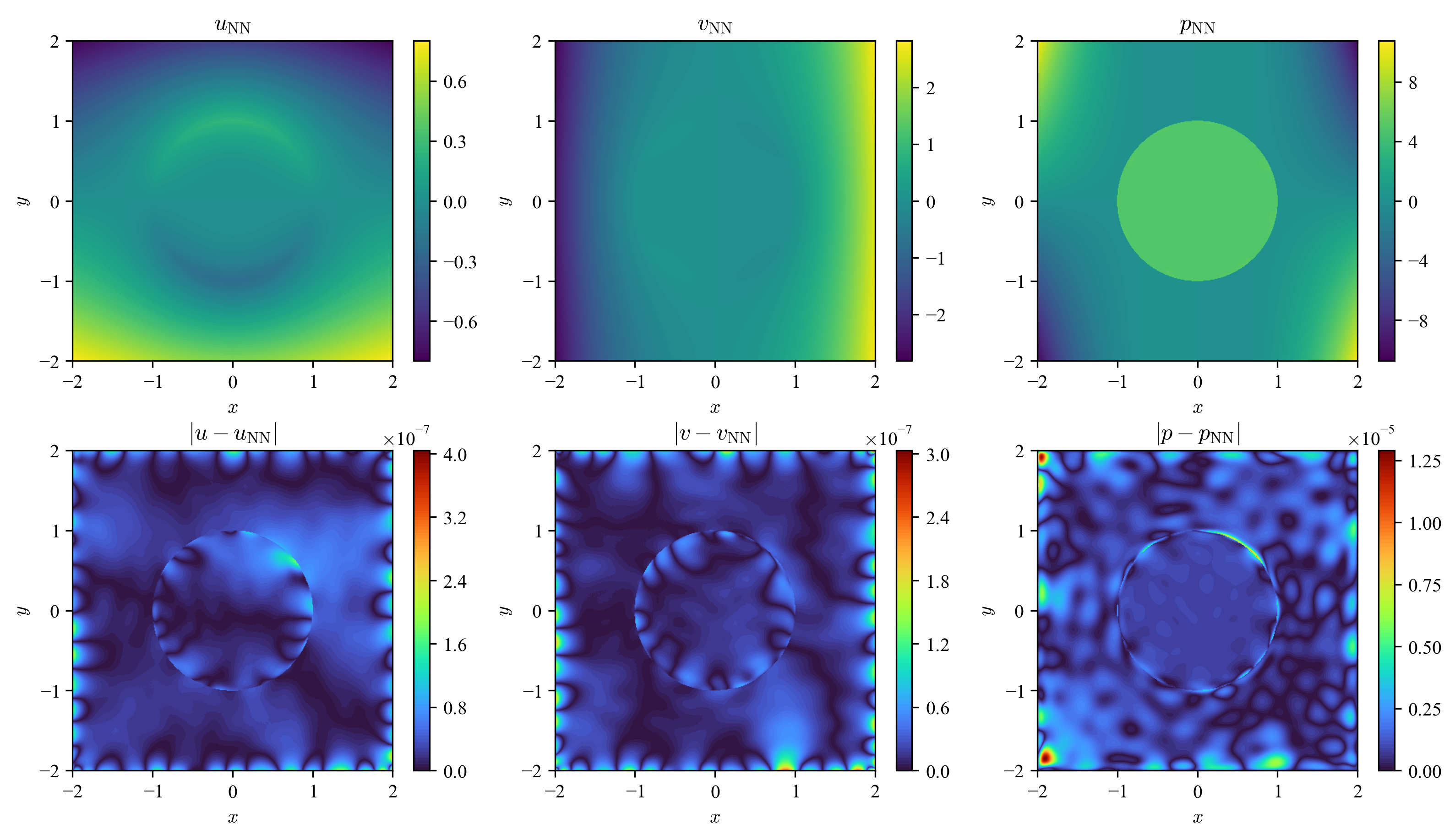}
     \vspace{-0.7cm}
    \caption{Numerical solutions (top row) and {pointwise absolute errors} (bottom row) produced by the LD-GTransNet method with $N=1600$ and $(\mu_1,\mu_2) = (1,10^3)$ for the 2D Stokes interface problem \eqref{eq:Stokes}.}
    \label{fig:Stokes}
\end{figure}

\begin{table}[ht!]
  \centering\small
  \caption{Performance comparison among LD-GTransNet, RFM \cite{Chi2024}, and DC-RFM \cite{SONG2026} for the 2D Stokes interface problem \eqref{eq:Stokes} under different viscosity contrasts.}
  \begin{tabular}{llccccc}
    \toprule
    $(\mu_1,\mu_2)$ & Method & \#Constraints & \#DOFs & $err_{L^2}(u)$ & $err_{L^2}(v)$
    & $err_{L^2}(p)$\\
    \midrule
    \multirow{3}{*}{$(1,10)$}
    & RFM & 84801 & 9600 & 6.92e-6 & 1.71e-6 & 3.33e-4 \\
    & DC-RFM & 10648 & 19200 & 1.79e-5 & 5.03e-6 & 6.73e-5 \\
    & LD-GTransNet & 12200 & 4800 & \textbf{9.90e-8} & \textbf{3.01e-8} & \textbf{4.39e-7} \\
    \midrule
    \multirow{3}{*}{$(1,10^3)$}
    & RFM & 84801 & 9600 & 4.49e-6 & 1.44e-6 & 3.31e-2 \\
    & DC-RFM & 10648 & 19200 & 1.20e-5 & 3.49e-6 & 2.91e-3 \\
    & LD-GTransNet & 12200 & 4800 & \textbf{6.23e-8} & \textbf{2.08e-8} & \textbf{1.69e-5}\\
    \bottomrule
  \end{tabular}
  \label{tab:comparison_Stokes}
\end{table}

First, Figure~\ref{fig:Stokes-errors} illustrates the relative \(L^{\infty}\) and \(L^2\) errors of the LD-GTransNet method versus the number of neurons \(N\) in the last hidden layer under different viscosity coefficients. It can be observed that the errors decrease significantly as \(N\) increases; when \(N=1600\), the relative \(L^2\) error of the velocity field \(\bu\) reaches the order of \(\mathcal{O}(10^{-7})\).
Next, Figure~\ref{fig:Stokes} presents the numerical solutions and the corresponding {pointwise absolute errors} of LD-GTransNet with \(N=1600\) and \((\mu_1,\mu_2) = (1,10^3)\).
Finally, we compare LD-GTransNet with RFM \cite{Chi2024} and DC-RFM \cite{SONG2026}. Table~\ref{tab:comparison_Stokes} reports the relative \(L^2\) errors under different viscosity contrasts. The results demonstrate that LD-GTransNet consistently achieves higher accuracy with a significantly smaller linear system across all test cases. {In particular, LD-GTransNet uses only 4,800 degrees of freedom--substantially fewer than RFM (9,600) and DC-RFM (19,200)--while improving accuracy by two to four orders of magnitude. The method remains robust across increasing viscosity contrasts (from $10$ to $10^3$), with only minor degradation in the pressure error, demonstrating that the decoupled architecture effectively handles large coefficient jumps without requiring additional tuning.}

\subsubsection{A 2D anisotropic diffusion-reaction {multiple-interface} problem}

\begin{figure}[ht!]
    \centering
    \includegraphics[width=1\linewidth]{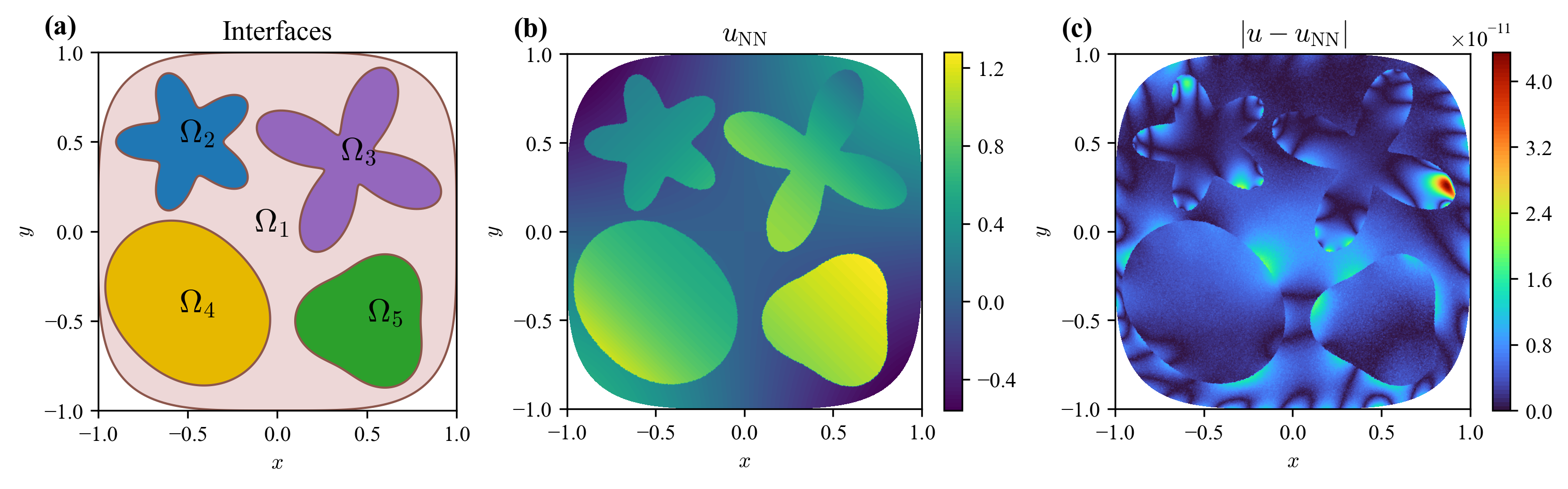}
     \vspace{-0.7cm}
    \caption{(a): Illustration of the partition of the domain into five subdomains by four interfaces. (b)-(c): Numerical solutions and {pointwise absolute errors} produced by the LD-GTransNet method with $N=1800$ and $(\kappa_1,\kappa_2,\kappa_3,\kappa_4,\kappa_5) = (1,10^{-1},10^{-2},10,10^2)$ for the 2D anisotropic {multiple-interface} problem \eqref{eq:anisotropic_multi_interfaces}.}
    \label{fig:anisotropic-multi-interface}
\end{figure}

In this subsection, we solve the following 2D anisotropic diffusion-reaction {multiple-interface} problem \cite{Lai2026} on the {superellipse} domain $\Omega = \{(x,y)\mid x^4 + y^4 \le 1\}$ with four interfaces. 
\begin{equation}\label{eq:anisotropic_multi_interfaces}
\begin{cases}
-\nabla \cdot (\Lambda\nabla u) + \mu u = f, & (x,y) \in \bigcup_{k=1}^5\Omega_k, \\
 \jump{u} = h_i, & (x,y) \in \Gamma_i,~i=1,2,3,4, \\
 \jump{\Lambda \nabla u \cdot \bn_i} = q_i, & (x,y) \in \Gamma_i,~i=1,2,3,4, \\
 u = g, & (x,y) \in \partial \Omega,
\end{cases}
\end{equation}
where the anisotropic diffusion tensor $\Lambda = \kappa_k\Lambda_0$ and the reaction coefficient $\mu = \kappa_k\mu_0$ for $(x,y) \in \Omega_k$ ($k=1,2,3,4,5$), with
\begin{equation}
\Lambda_0 = \begin{bmatrix}
    (x+y)^2 + 1 & y^2 - x^2 \\
    y^2 - x^2 & (x-y)^2 + 1
\end{bmatrix},
\quad
\mu_0 = \exp(x-y).
\end{equation}
As illustrated in Figure~\ref{fig:anisotropic-multi-interface} (a), the four interfaces are defined as
\begin{align*}
    \Gamma_1 &= \{(x,y)\mid (x+0.5)^2 + (y-0.5)^2 = (0.3-0.1\cos(5\theta))^2\}, \\
    \Gamma_2 &= \{(x,y)\mid (x-0.4)^2 + (y-0.4)^2 = (0.35-0.2\sin(4\theta))^2\}, \\
    \Gamma_3 &= \{(x,y)\mid (x+0.5)^2 + (y+0.4)^2 = (0.45-0.05\sin(2\theta))^2\}, \\
    \Gamma_4 &= \{(x,y)\mid (x-0.5)^2 + (y+0.5)^2 = (0.35-0.05\cos(3\theta))^2\},
\end{align*}
where $\theta \in [0,2\pi)$. The exact solution is given by
\[
u(x,y) = \begin{cases}
    \sin x\sin y, & (x,y)\in\Omega_1, \\
    \exp(x-y), & (x,y)\in\Omega_2, \\
    \cos(x+y), & (x,y)\in\Omega_3, \\
    0.5\cosh(x+y), & (x,y)\in\Omega_4,\\
    \ln(x+y+3), & (x,y)\in\Omega_5.
\end{cases}
\]
The source term $f$, boundary condition $g$, and jump conditions $h_i, q_i$ ($i=1,2,3,4$) are again determined {from the exact solution}.

We set \(\xi\) as
\[
\xi =
\begin{cases}
1, & (x,y) \in \Omega_1, \\
0.5, & (x,y) \in \Omega_2, \\
0, & (x,y) \in \Omega_3, \\
-0.5, & (x,y) \in \Omega_4, \\
-1, & (x,y) \in \Omega_5,
\end{cases}
\]
and $(\mathbf{x}_c,R,\xi_c,R_{\xi}) = ((0,0),1.5,0,1.1)$.
The number of neurons in the first hidden layer of LD-GTransNet is fixed at 1800.
The collocation points are deployed as follows: 1684 points in $\Omega$, 200 points on $\partial\Omega$, and 100 points on each interface $\Gamma_1, \Gamma_2, \Gamma_3, \Gamma_4$, respectively.

\begin{figure}[ht!]
    \centering
    \includegraphics[width=0.8\linewidth]{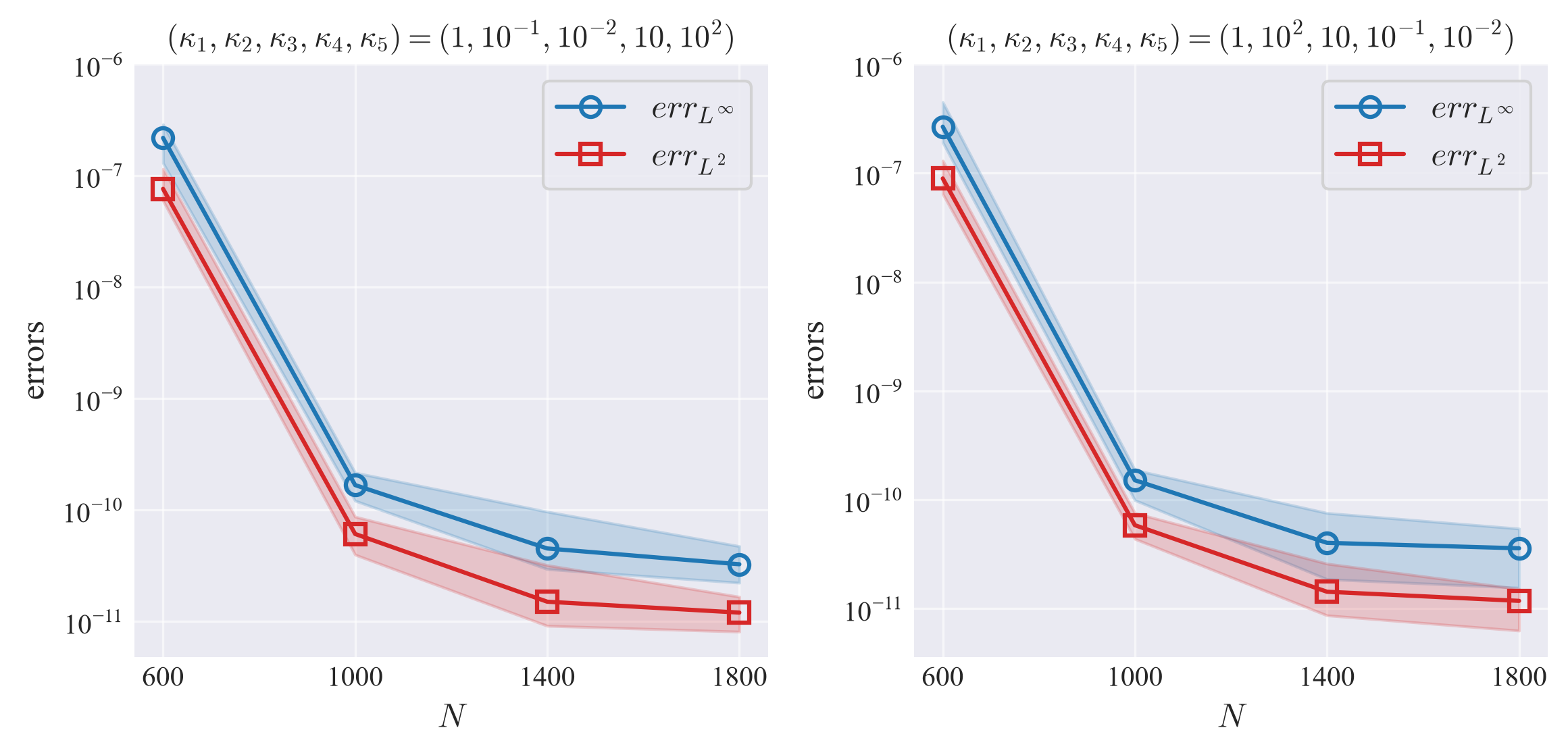}
     \vspace{-0.2cm}
    \caption{Relative $L^{\infty}$ and $L^2$ errors of the numerical solutions obtained by the LD-GTransNet method with respect to the number of neurons $N$ in the last hidden layer for the 2D anisotropic diffusion-reaction {multiple-interface} problem \eqref{eq:anisotropic_multi_interfaces} under different combinations of diffusion coefficients.}
    \label{fig:anisotropic-multi-interface-errors}
\end{figure}

\begin{table}[ht!]
  \centering\small
  \caption{Performance comparison between LD-GTransNet and DCSNN \cite{Hu2022} for the 2D anisotropic diffusion-reaction {multiple-interface} problem \eqref{eq:anisotropic_multi_interfaces} under different combinations of diffusion coefficients.}
  \begin{tabular}{llccc}
    \toprule
    $(\kappa_1, \kappa_2, \kappa_3, \kappa_4,\kappa_5)$ & Method & \#Constraints & \#DOFs & $err_{L^2}(u)$ \\
    \midrule
    \multirow{2}{*}{$(1,10^{-1}, 10^{-2}, 10, 10^2)$}
    & DCSNN & 2684 & - & 3.41e-7 \\
    & LD-GTransNet & 2684 & 1600 & \textbf{8.30e-12} \\
    \midrule
    \multirow{2}{*}{$(1,10^{2}, 10, 10^{-1}, 10^{-2})$}
    & DCSNN & 2684 & - & 4.56e-7 \\
    & LD-GTransNet & 2684 & 1600 & \textbf{8.67e-12} \\
    \bottomrule
  \end{tabular}
  \label{tab:comparison_anisotropic_multi_interface}
\end{table}

Similarly, Figure~\ref{fig:anisotropic-multi-interface-errors} presents the variations of the relative $L^{\infty}$ and $L^2$ errors of LD-GTransNet with respect to the number of neurons $N$ in the last hidden layer under different combinations of anisotropic diffusion coefficients.
It can be observed that the relative $L^{\infty}$ and $L^2$ errors of LD-GTransNet reach $\mathcal{O}(10^{-11})$ and $\mathcal{O}(10^{-12})$ at $N=1800$, respectively. Table~\ref{tab:comparison_anisotropic_multi_interface} compares the performance of LD-GTransNet with that of DCSNN.
{DCSNN~\cite{Hu2022} is a fully parameterized neural network whose trainable parameters are optimized via the Levenberg-Marquardt (LM) algorithm and the nonlinear optimization incurs a significantly higher computational cost than the linear least-squares solve of LD-GTransNet.} The results demonstrate that LD-GTransNet substantially outperforms DCSNN in terms of numerical accuracy. Finally, Figure~\ref{fig:anisotropic-multi-interface} (b)-(c) visualize the numerical solutions and the corresponding {pointwise absolute errors} of LD-GTransNet with $N=1800$.

\subsubsection{A 3D elasticity interface problem}\label{subsec:elasticity}

Finally, we solve the following 3D elasticity interface problem \cite{HO_NIPFEM} on the domain $\Omega = (0,1)^3$:
\begin{equation}\label{eq:elasticity}
\begin{cases}
-\nabla \cdot \boldsymbol{\sigma}(\bu) = \mathbf{f}, & (x,y,z) \in \Omega_1 \cup \Omega_2, \\
\quad\quad\ \, \jump{\bu} = \mathbf{h}, & (x,y,z) \in \Gamma, \\
\quad \jump{\boldsymbol{\sigma}(\bu)\bn} = \mathbf{q}, & (x,y,z) \in \Gamma, \\
\quad\quad\quad\ \bu = \mathbf{g}, & (x,y,z) \in \partial \Omega,
\end{cases}
\end{equation}
where the displacement vector \( \bu = (u, v, w)^\top \), the stress tensor \( \boldsymbol{\sigma}(\bu) = \nu (\nabla \cdot \bu) \mathbf{I} + 2\mu \boldsymbol{\varepsilon}(\bu) \), the strain tensor \( \boldsymbol{\varepsilon}(\bu) = \left(\nabla \bu + \nabla \bu^\top\right)/2 \), and the Lam\'e parameters \( \nu(x,y,z) = \nu_k \) and \( \mu(x,y,z) = \mu_k \) for \( (x,y,z) \in \Omega_k \) with \(k = 1,2\).
The interface is defined as
\begin{equation}
\Gamma = \left\{(x,y,z)\mid\frac{(x-0.5)^2}{0.15^2} + \frac{(y-0.5)^2}{0.2^2} + \frac{(z-0.5)^2}{0.25^2} = 1 \right\}.
\end{equation}
Let $E(x,y,z)=(x-0.5)^2/0.15^2+(y-0.5)^2/0.2^2+(z-0.5)^2/0.25^2$. The subdomains are $\Omega_1=\{(x,y,z)\in\Omega\mid E<1\}$ and $\Omega_2=\{(x,y,z)\in\Omega\mid E>1\}$, with $(\nu,\mu)=(\nu_1,\mu_1)$ inside the ellipsoid and $(\nu,\mu)=(\nu_2,\mu_2)$ outside. The exact solutions are given by
\[
u(x,y,z) = 
\begin{cases}
-\cos\left(x^2\right)\mathrm{e}^{-y^2}\sin(2\pi z), & (x,y,z) \in \Omega_1, \\
-\sin\left(x^2\right)\mathrm{e}^{y^2}\cos(2\pi z), & (x,y,z) \in \Omega_2,
\end{cases}
\]
\[
v(x,y,z) = 
\begin{cases}
-\cos\left(y^2\right)\mathrm{e}^{-x^2}\sin(2\pi z), & (x,y,z) \in \Omega_1, \\
-\sin\left(y^2\right)\mathrm{e}^{x^2}\cos(2\pi z), & (x,y,z) \in \Omega_2,
\end{cases}
\]
and
\[
w(x,y,z) = 
\begin{cases}
\cos\left(y^2\right)\mathrm{e}^{-z^2}\sin(2\pi x), & (x,y,z) \in \Omega_1, \\
\sin\left(y^2\right)\mathrm{e}^{z^2}\cos(2\pi x), & (x,y,z) \in \Omega_2.
\end{cases}
\]
The source term $\mathbf{f}$, boundary condition $\mathbf{g}$, and jump conditions $\mathbf{h}, \mathbf{q}$ are computed from the exact solutions.

\begin{figure}[ht!]
    \centering
    \includegraphics[width=0.8\linewidth]{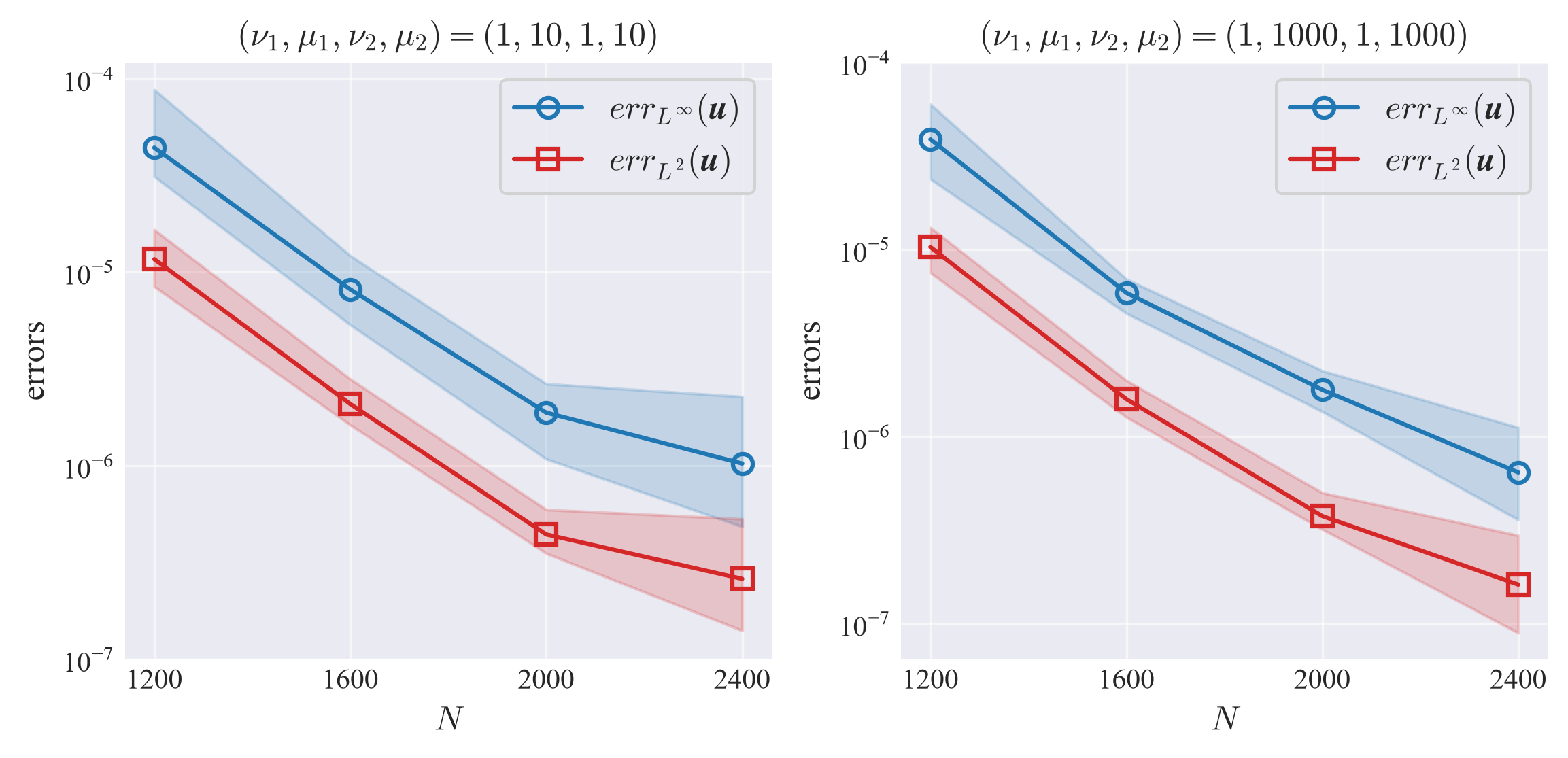}
     \vspace{-0.2cm}
    \caption{Relative $L^{\infty}$ and $L^2$ errors of the numerical solutions obtained by the LD-GTransNet method with respect to the number of neurons $N$ in the last hidden layer for the 3D elasticity interface problem \eqref{eq:elasticity} with different Lam\'e parameters.}
    \label{fig:elasticity-errors}
\end{figure}

\begin{table}[ht!]
  \centering\small
  \caption{Performance comparison between LD-GTransNet and high-order NIPFEM \cite{HO_NIPFEM} for the 3D elasticity interface problem with high-contrast Lam\'e parameters $(\nu_1,\mu_1,\nu_2,\mu_2) = (1,1000,1,1000)$.}
  \begin{tabular}{lccc}
    \toprule
    Method & \#Constraints & \#DOFs & $err_{L^2}(\bu)$ \\
    \midrule
    \multirow{3}{*}{NIPFEM} & 23500 & 23500 & 2.96e-3 \\
           & 176694 & 176694 & 1.38e-4 \\
           & 1369446 & 1369446 & 6.96e-6 \\
    LD-GTransNet & 50580 & 6000 & \textbf{2.60e-7} \\
    \bottomrule
  \end{tabular}
  \label{tab:comparison_elasticity}
\end{table}

\begin{figure}[ht!]
    \centering
    \includegraphics[width=1\linewidth]{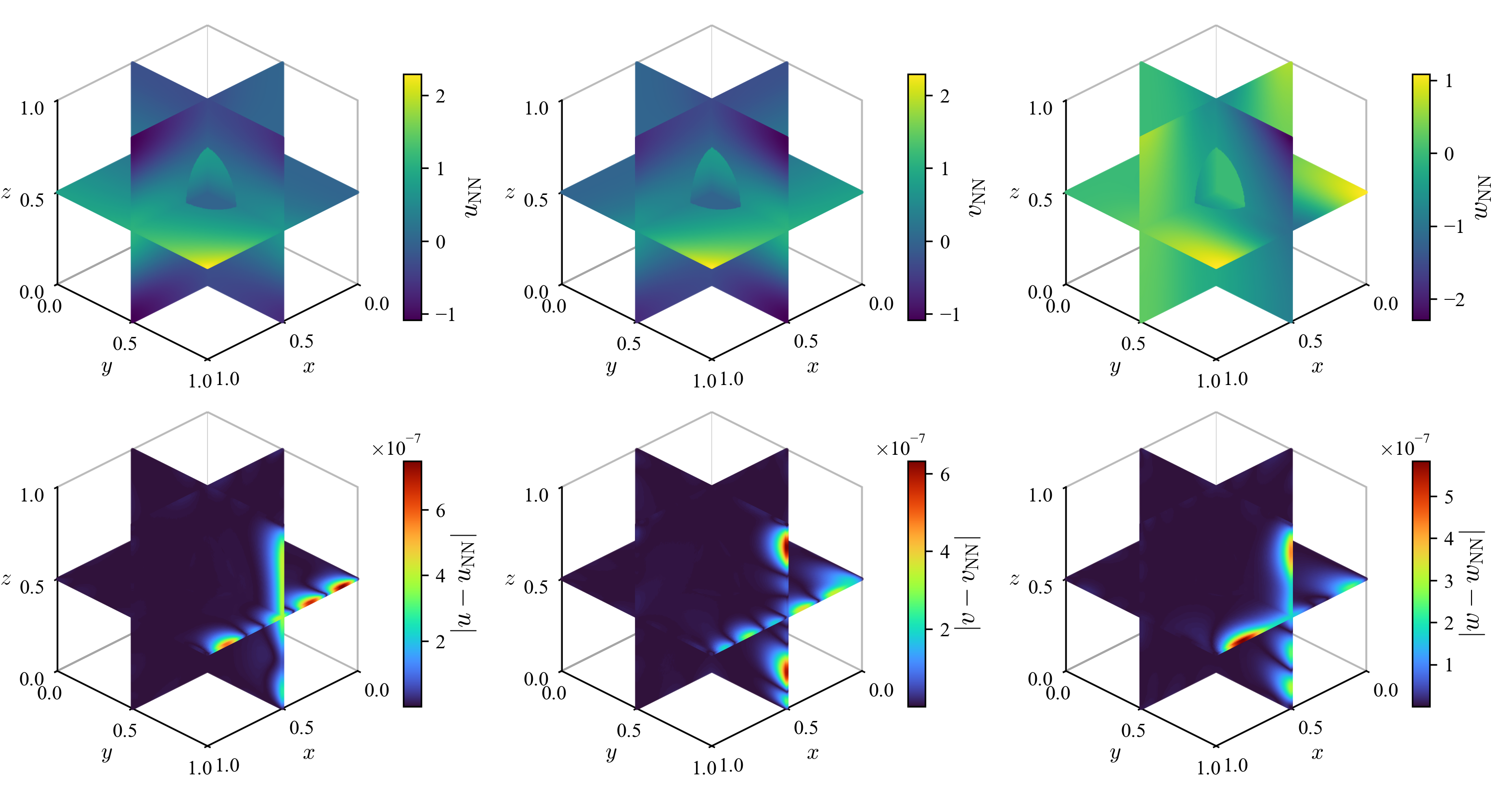}
     \vspace{-0.6cm}
    \caption{Numerical solution (top row) and pointwise absolute error (bottom row) produced by the LD-GTransNet method with $N=2400$ and $(\nu_1,\mu_1,\nu_2,\mu_2) = (1,1000,1,1000)$ for the 3D elasticity interface problem \eqref{eq:elasticity} at $x=0.5$, $y=0.5$ and $z=0.5$.}
    \label{fig:elasticity}
\end{figure}

We set \(\xi\) as
\[
\xi =
\begin{cases}
1, & (x,y,z) \in \Omega_1, \\
-1, & (x,y,z) \in \Omega_2,
\end{cases}
\]
and the parameters are set to $\mathbf{x}_c=(0.5,0.5,0.5),R=1$ for GTransNet and $\xi_c=0,R_{\xi}=1.1$ for TransNet.   
The number of neurons in the first hidden layer of LD-GTransNet is fixed at 3000.
The collocation points are deployed as follows: 12000 points in $\Omega$, 1350 points on $\partial\Omega$, and 1800 points on the interface $\Gamma$, respectively.

Consistent with the 2D numerical examples presented above, Figure~\ref{fig:elasticity-errors} first illustrates the relative $L^{\infty}$ and $L^2$ errors of the LD-GTransNet method with respect to the number of neurons in the last hidden layer, and the results exhibit similar convergence behavior to the 2D cases.
Meanwhile, Table~\ref{tab:comparison_elasticity} presents the performance comparison between the LD-GTransNet method and the high-order NIPFEM method for the 3D elasticity interface problem.
The LD-GTransNet method demonstrates superior performance in terms of the number of constraints, degrees of freedom, and $L^2$ error.
Finally, Figure~\ref{fig:elasticity} displays the numerical solution and pointwise absolute error of the LD-GTransNet method for the 3D elasticity interface problem.

\subsection{Tests for moving interface problems}
In this subsection, we use the proposed LD-GTransNet for solving some time-dependent moving interface problems. 

\subsubsection{A 2D moving interface problem}
We first consider a 2D moving interface problem on the domain $\Omega = (-1,1)^2$, featuring {an evolving circular interface} {defined} by:
\begin{equation}
\Gamma(t) = \{(x,y)\mid x^2 + y^2 = 0.5+0.3t\}.
\end{equation}
The subdomains are $\Omega_1(t)=\{(x,y)\in\Omega\mid x^2+y^2<0.5+0.3t\}$ and $\Omega_2(t)=\{(x,y)\in\Omega\mid x^2+y^2>0.5+0.3t\}$. The diffusion coefficient is $\kappa_1$ inside the circle and $\kappa_2$ outside; $\bn$ points from $\Omega_1(t)$ to $\Omega_2(t)$. The problem is governed by the following parabolic equation:
\begin{equation}\label{eq:moving_interface}
\begin{cases}
\ds \frac{\partial u}{\partial t}-\nabla \cdot (\kappa \nabla u) = f, & (x,y) \in \Omega_1(t)\cup\Omega_2(t),~t \in(0,T], \\
\jump{u} = h, & (x,y) \in \Gamma(t),~t \in (0,T], \\
\jump{\kappa \nabla u \cdot \bn} = q, & (x,y) \in \Gamma(t),~t\in(0,T], \\
u = g, & (x,y) \in \partial \Omega,~t\in (0,T], \\
u = u_0, & (x,y)\in\Omega,~t = 0.
\end{cases}
\end{equation}
where the diffusion coefficient $\kappa$ takes the constant value $\kappa_k$ in each subdomain $\Omega_k(t)$ ($k=1,2$), and the terminal time is set to $T=1$.
We set the exact solution as:
\begin{equation*}
    u(x,y,t) = \begin{cases}
       \cos(1 - x^2 - y^2 - t) & (x,y,t) \in\Omega_1(t),\\ 
       \log(1+x^2+y^2)\sin t,& (x,y,t) \in \Omega_2(t).
    \end{cases}
\end{equation*}
The source term $f$, boundary condition $g$, initial condition $u_0$ and jump conditions $h, q$ are again determined {from the exact solution}. We set \(\xi\) as
\[
\xi =
\begin{cases}
1, & (x,y,t) \in \Omega_1(t), \\
-1, & (x,y,t) \in \Omega_2(t),
\end{cases}
\]
and the parameters are set to $(x_c,y_c,t_c)=(0,0,0.5)$, $R=2$ for GTransNet and $\xi_c=0$, $R_{\xi}=1.1$ for TransNet.  The number of neurons in the first hidden layer of LD-GTransNet is fixed at 3000.
Collocation points are sampled once in the entire space-time domain $\Omega\times(0,T]$ using Latin hypercube sampling and are not regenerated at individual time steps; the time variable $t$ is treated as an additional spatial coordinate, so the moving interface becomes a fixed hypersurface in the lifted space. The collocation points are deployed as follows: 4000 points in $\Omega(t)$, 1600 points on $\partial\Omega(t)$, 400 points on $\Omega(0)$ and 1800 points on the interface $\Gamma(t)$, respectively.

\begin{figure}[ht!]
    \centering
    \includegraphics[width=0.9\linewidth]{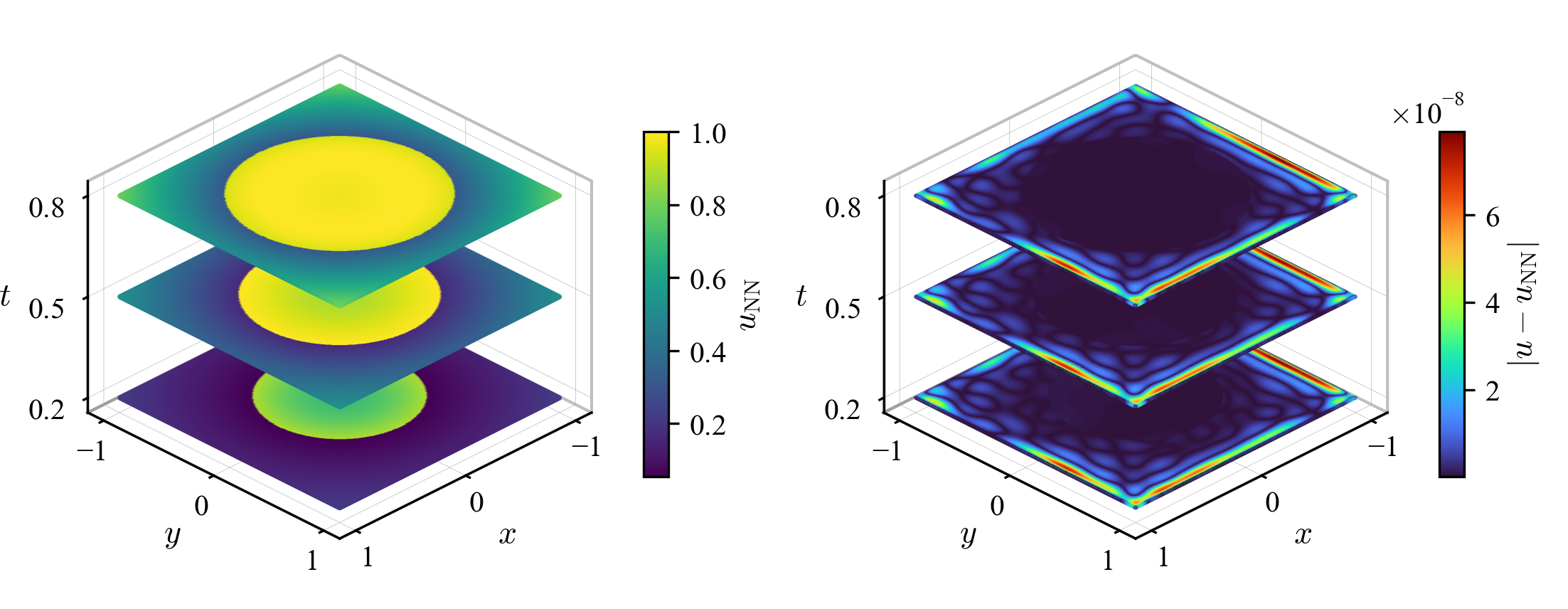}
     \vspace{-0.2cm}
    \caption{Numerical solutions and {pointwise absolute errors} produced by the LD-GTransNet method with $N=3000$ and $(\kappa_1,\kappa_2) = (1,10^{-6})$ for the moving interface problem \eqref{eq:moving_interface} at $t=0.2,0.5$ and $0.8$.}
    \label{fig:evolving_sphere}
\end{figure}

\begin{table}[ht!]
  \centering\small
  \caption{Performance comparison between LD-GTransNet and INN \cite{Wu2022} for the moving interface problem \eqref{eq:moving_interface} under different combinations of diffusion coefficients.}
  \begin{tabular}{llccc}
    \toprule
    $(\kappa_1, \kappa_2)$ & Method & \#Constraints & \#DOFs & $err_{L^2}$ \\
    \midrule
    \multirow{2}{*}{$(1,10^{-6})$}
    & INN & 9600 & - & 8.43e-4 \\
    & LD-GTransNet & 9600 & 3000 & \textbf{1.28e-8} \\
    \midrule
    \multirow{2}{*}{$(1,10^6)$}
    & INN & 9600 & - & 6.56e-4 \\
    & LD-GTransNet & 9600 & 3000 & \textbf{2.77e-9} \\
    \bottomrule
  \end{tabular}
  \label{tab:comparison_moving_interface}
\end{table}

Figure~\ref{fig:evolving_sphere} presents the numerical solutions and {pointwise absolute errors} obtained by the LD-GTransNet method for the moving interface problem. Table~\ref{tab:comparison_moving_interface} compares the performance of the LD-GTransNet method with that of the Interface Neural Networks method \cite{Wu2022} on the moving interface problem. {Notably, LD-GTransNet eliminates nonlinear optimization entirely---INN requires gradient-based training---while achieving four to five orders of magnitude higher accuracy using the same number of constraints.}

\subsubsection{A 2D {two-phase} Stefan problem}
Consider a 2D {two-phase} Stefan problem \eqref{eq:stefan} on the domain $\Omega = [0,2]\times[0,1]$ with terminal time $T=1$, where
$\Omega_1(t)=\{(x,y)\in(0,2)\times(0,1)\mid x<s(y,t)\}$ and $\Omega_2(t)=\{(x,y)\in(0,2)\times(0,1)\mid x>s(y,t)\}$. The diffusion coefficient is $\kappa_1=2$ in $\Omega_1(t)$ and $\kappa_2=1$ in $\Omega_2(t)$. The interface \(\Gamma(t)=\{(x,y)\in\Omega\mid x=s(y,t)\}\) is determined by
\begin{equation}
\begin{cases}
\ds \frac{\partial s(y,t)}{\partial t} = -2\left(\frac{\partial u_1(s,y,t)}{\partial x}
- \frac{\partial u_1(s,y,t)}{\partial y} \frac{\partial s(y,t)}{\partial y}\right) + \\
\qquad\qquad\qquad\ds\left(\frac{\partial u_2(s,y,t)}{\partial x}
- \frac{\partial u_2(s,y,t)}{\partial y} \frac{\partial s(y,t)}{\partial y}\right), \\
s(y,0) = \frac12.
\end{cases}
\end{equation}
The exact solution is given by
\begin{equation*}
    u(x,y,t) = \begin{cases}
        \ds 2\exp(\frac{2t-2x+1}{4}) - 2,& (x,y,t) \in \Omega_1(t),\\
        \ds \exp(\frac{2t-2x+1}{2}) - 1, & (x,y,t) \in \Omega_2(t).
    \end{cases}
\end{equation*}

We set \(\xi\) as
\[
\xi =
\begin{cases}
1, & (x,y,t) \in \Omega_1(t), \\
-1, & (x,y,t) \in \Omega_2(t),
\end{cases}
\]
and $(\mathbf{x}_c,R,\xi_c,R_{\xi}) = ((0.5,1,0.5),1.5,0,1.1)$.
The number of neurons in the first hidden layer of LD-GTransNet is fixed at 2000, and the 
number of neurons in the hidden layer of TransNet for $s(y,t)$ is fixed at 1000.
The collocation points are deployed as follows: 2000 points in $\Omega(t)$, 400 points on $\partial\Omega(t)$, 100 points on $\Omega(0)$ and 400 points on the interface $\Gamma(t)$, respectively. 

\begin{figure}[ht!]
    \centering
    \includegraphics[width=0.9\linewidth]{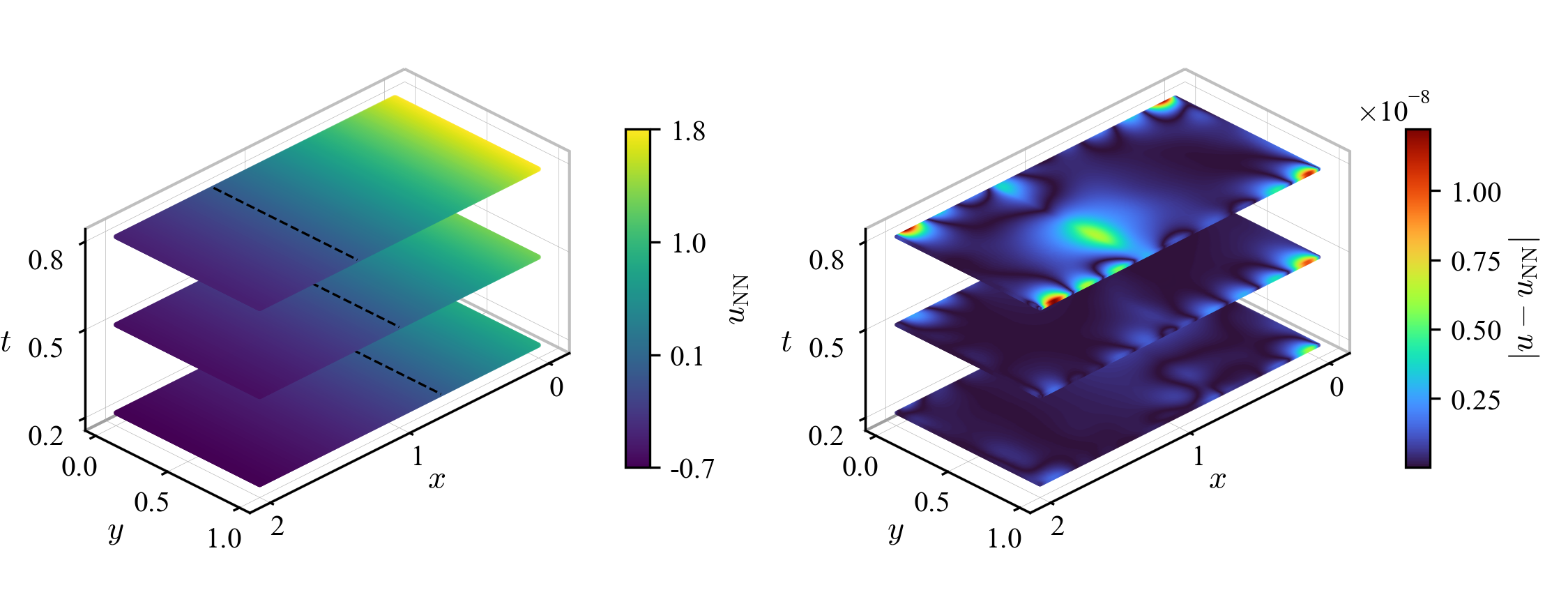}
     \vspace{-0.3cm}
    \caption{Numerical solutions and {pointwise absolute errors} produced by the LD-GTransNet method with $N=2000$ for the two-phase Stefan problem \eqref{eq:stefan} at $t=0.2,0.5$ and $0.8$.}
    \label{fig:stefan}
\end{figure}

\begin{table}[ht!]
  \centering\small
  \caption{Performance comparison between LD-GTransNet and PIELM \cite{PIELM_Stefan} for the {two-phase} Stefan problem \eqref{eq:stefan}.}
  \begin{tabular}{lccc}
    \toprule
    Method & \#Constraints & \#DOFs & $err_{L^2}$ \\
    \midrule
     PIELM & 3300 & 2000 & 2.74e-6 \\
    LD-GTransNet & 3300 & 2000 & \textbf{1.81e-9} \\
    \bottomrule
  \end{tabular}
  \label{tab:comparison_stefan}
\end{table}

Figure~\ref{fig:stefan} displays the numerical solutions and {pointwise absolute errors} of the LD-GTransNet method for the Stefan problem.
Furthermore, Table~\ref{tab:comparison_stefan} presents the performance comparison between the LD-GTransNet method and the PIELM method \cite{PIELM_Stefan} for the Stefan problem. {With the same number of constraints and degrees of freedom, LD-GTransNet achieves an $L^2$ error three orders of magnitude lower than PIELM, demonstrating the effectiveness of the lift-and-decoupling strategy for coupled nonlinear moving-interface problems.}

\section{Conclusion}\label{sec:5}
We proposed {a novel lift-and-decoupling} generalized transferable neural network for solving PDEs with interface discontinuities. By decoupling physical and auxiliary variables, it provides a unified global representation of piecewise smooth solutions without domain decomposition, avoiding the redundancy issue.
The framework reduces to a linear least-squares problem with fixed hidden-layer parameters, eliminating nonlinear optimization and enabling efficient, robust, and easy to implement solvers. It is broadly applicable to elliptic equations, Stokes flow, anisotropic diffusion–reaction systems, elasticity problems, and moving interfaces.
Numerical experiments in two and three dimensions show that LD-GTransNet achieves high-order accuracy, while using significantly smaller linear systems than existing methods. Its advantages are especially pronounced for problems with high-contrast coefficients or complex, evolving interfaces.
Future work will focus on theoretical analysis, extensions to high-dimensional and nonlinear systems, and integration with adaptive sampling or multilevel techniques to enhance scalability and efficiency.

\section*{Acknowledgements}
L. Ju's research is partially supported
by U.S. Department of Energy grant DE-SC0025527. X. Zhang's research is partially supported by Natural Science Foundation of China grant 12571467.

\bibliographystyle{plain}
\bibliography{ref}
\end{document}

%% file: preamble.tex
\newcommand{\ba}{\mathbf{a}}

\newcommand{\bn}{\mathbf{n}}

\newcommand{\br}{\mathbf{r}}

\newcommand{\bu}{\mathbf{u}}

\newcommand{\bw}{\mathbf{w}}
\newcommand{\bx}{\mathbf{x}}

\newcommand{\bA}{\mathbf{A}}

\newcommand{\bW}{\mathbf{W}}

\newcommand{\cB}{\mathcal{B}}

\newcommand{\cD}{\mathcal{D}}

\newcommand{\cL}{\mathcal{L}}

\newcommand{\valpha}{{\boldsymbol{\alpha}}}

\newcommand{\vvarphi}{\boldsymbol{\varphi}}
\newcommand{\vpsi}{\boldsymbol{\psi}}

\newcommand{\R}{\mathbb R}
\newcommand{\ds}{\displaystyle}

\newcommand{\jump}[1]{\ensuremath{[\![ #1 ]\!]}}

\definecolor{blue}{rgb}{0.0,0.0,1.0}
\definecolor{red}{rgb}{1.0,0.0,0.0}
\definecolor{purple}{rgb}{0.75, 0.0, 1.0}